\documentclass[12pt,reqno]{amsart}%{article}
\usepackage{pb-diagram} 
\usepackage{amsmath,amssymb}
\usepackage{enumerate}
\usepackage{amsthm}
\usepackage{hyperref}
\usepackage{pb-diagram}

\newtheorem{theorem}{Theorem}
\newtheorem{proposition}[theorem]{Proposition}

\newtheorem{remark}{Remark}
\newtheorem{lemma}[theorem]{Lemma}
\newfont{\bb}{msbm10 at 12pt}

\def\Ga{{\Gamma}}

\def\phi{{\varphi}}

\def\M{\mathrm{M}}
\def\T*M{\mathrm{T}^*\mathrm{M}}

\def\gst{{g_{\rm st}}}

\def\hs{\mathbb{S}^n_+}

\def\V{\mathrm{V}}

\def\ov{\overline}

\let\pa\partial     
\let\na\nabla     
     
\DeclareMathAlphabet{\doba}{U}{msb}{m}{n}

\def\T{\mathrm{T}}
\def\Vol{{\mathop{\rm Vol}}}

\def\Vol{{\mathop{\rm Vol}}}

\def\bhl{\lambda_{\min}(M,[g],\sigma)}

\let\<\langle     
\let\>\rangle

\let\<\langle     
\let\>\rangle

\long\def\komment#1{} 

\begin{document}

%%%%%%%%%%%%%%%%%%%%%%%%%%%%%%%%%%%%%%%%%%%%%%%%%%%%%%%%%%%%%%%%%%%%%%%%%%%%%%%%%%%%%%%%%%%%%%%%%%%%%%%%%     

\title{A Sobolev-like Inequality for the Dirac operator}     
\author{Simon Raulot}\thanks{The author is supported by the Swiss SNF grant $20$-$118014/1$}
\address{Universit\'e de Neuch\^atel\\Institut de Math\'ematiques,  
Rue Emile-Argand 11 \\2007 Neuch\^atel\\ Suisse}
\email{simon.raulot@unine.ch}  
\date{\today}
\keywords{Dirac operator, Sobolev inequality, Conformal Geometry, Nonlinear elliptic equations}
\subjclass{}

\maketitle

\begin{abstract}
In this article, we prove a Sobolev-like inequality for the Dirac operator on closed compact Riemannian spin manifolds with a nearly optimal Sobolev constant. As an application, we give a criterion for the existence of solutions to a nonlinear equation with critical Sobolev exponent involving the Dirac operator. We finally specify a case where this equation can be solved.
\end{abstract}

%%%%%%%%%%%%%%%%%%%%%%%%%%%%%%%%%%%%%%%%%%%%%%%%%%%%%%%%%%%%%%%%%%%%%%%%%%%%%%%%%%%%%%%%%%%%%%%%%%%%%%%%

\section{Introduction}

%%%%%%%%%%%%%%%%%%%%%%%%%%%%%%%%%%%%%%%%%%%%%%%%%%%%%%%%%%%%%%%%%%%%%%%%%%%%%%%%%%%%%%%%%%%%%%%%%%%%%%%%

Let $(M^n,g)$ be a compact Riemannian manifold of dimension $n\geq 3$. The Sobolev embedding theorem asserts that the Sobolev space $H^2_1$ of functions $u\in L^2$ such that $\nabla u\in L^2$ embeds continuously in the Lebesgue space $L^N$ (with $N=\frac{2n}{n-2}$). In other words, there exists two constants $A,B>0$ such that, for all $u\in H^2_1$, we have:
$$
\Big(\int_M|u|^Ndv(g)\Big)^{\frac{2}{N}}\leq A\int_{M}|\nabla u|^2dv(g)+B\int_Mu^2dv(g).\eqno S(A,B)
$$ 

Considerable work has been devoted to the analysis of sharp Sobolev-type inequalities, very often in connection with concrete problems from geometry. One of these concerns the best constant in $S(A,B)$ defined by:
\begin{eqnarray*}
A_2(M):=\inf\mathcal{A}_2(M)
\end{eqnarray*}

where: 
\begin{eqnarray*}
\mathcal{A}_2(M):=\{A>0\,/\,\exists B>0\,\textrm{such that } S(A,B) \textrm{ holds for all } u\in C^\infty(M)\}.
\end{eqnarray*}

From $S(A,B)$ and by definition of $A_2$, we easily get that:
\begin{enumerate}
\item $A_2(M)\geq K(n,2)^2$

\item for any $\varepsilon>0$ there exists $B_\varepsilon>0$ such that inequality $S(A_2(M)+\varepsilon,B_\varepsilon)$ holds.
\end{enumerate}

Here $K(n,2)^2$ denotes the best constant of the corresponding Sobolev embedding theorem in the Euclidean space given by (see \cite{aubin} and \cite{Talenti}):
\begin{eqnarray*}
K(n,2)^2=\frac{4}{n(n-2)\omega_n^{\frac{2}{n}}}
\end{eqnarray*}
 
where $\omega_n$ stands for the volume of the standard $n$-dimensional sphere. In fact, Aubin \cite{aubin} showed that $A_2(M)=K(n,2)^2$ and conjectured that $S(A,B)$ should hold for $A=K(n,2)^2$, that is $\mathcal{A}_2(M)$ is closed. The proof of this conjecture by Hebey and Vaugon (see \cite{hebeyvaugon} and \cite{hebeyvaugon2}) gave rise to various interesting problems dealing with best constants in Riemannian Geometry. One of those given in \cite{hebmatz}, is the problem of prescribed critical functions which study the existence of functions for which $S(A_2(M),B_0)$ is an equality (here $B_0>0$ denotes the infimum on $B>0$ such that $S(A_2(M),B_0)$ holds). For more details and related topics, we refer to \cite{heb}.

Recall that one of the first geometric applications of the best constant problem has been discovered by Aubin \cite{aubin:76} regarding the Yamabe problem. This famous problem of Riemannian geometry can be stated as follows: given a compact Riemannian manifold $(M^n,g)$ of dimension $n\geq 3$, can one find a metric conformal to $g$ such that its scalar curvature is constant? This problem has a long and fruitful history and it has been completely solved in several steps by Yamabe \cite{yamabe:60}, Tr\"udinger \cite{trudinger:68}, Aubin \cite{aubin:76} and finally Schoen \cite{schoen:84} using the Positive Mass Theorem coming from General Relativity (see also \cite{lee.parker:87} for a complete review). The Yamabe problem is in fact equivalent to find a smooth positive solution $u\in C^\infty(M)$ to a nonlinear elliptic equation:
\begin{eqnarray}\label{yamequ}
L_g u:=4\frac{n-1}{n-2}\Delta_g u+R_gu=\lambda u^{N-1},
\end{eqnarray}  

where $L_g$ is known as the conformal Laplacian (or the Yamabe operator), $\Delta_g$ (resp. $R_g$) denotes the standard Laplacian acting on functions (resp. the scalar curvature) with respect to the Riemannian metric $g$ and $\lambda\in\mathbb{R}$ is a constant. Indeed, if such a function exists then the metric $\ov{g}=u^{N-2}g$ is conformal to $g$ and satisfies $R_{\ov{g}}=\lambda$. A standard variational approach cannot allow to conclude because of the lack of compactness in the Sobolev embedding theorem involved in this method. However, Aubin \cite{aubin:76} proved that if:
\begin{eqnarray}\label{ys}
Y(M,[g])=\inf_{f\neq 0} I(f)<Y(\mathbb{S}^n,[g_{st}])=4\frac{n-1}{n-2}K(n,2)^{-2}
\end{eqnarray} 

holds, where $I$ denotes the functional defined by:
\begin{eqnarray*}
I(f)=\frac{4\frac{n-1}{n-2}\int_M|\nabla f|^2dv(g)+\int_MR_gf^2dv(g)}{\Big(\int_M|f|^Ndv(g)\Big)^{\frac{2}{N}}},
\end{eqnarray*}

Equation (\ref{yamequ}) admits a positive smooth solution. This condition points out the tight relation between the Yamabe problem and the best constant involved in the Sobolev inequality. Moreover, it is sharp in the sense that for all compact Riemannian manifolds $(M^n,g)$, the following inequality holds (see \cite{aubin:76}):
\begin{eqnarray}\label{aubinin}
Y(M,[g])\leq 4\frac{n-1}{n-2}K(n,2)^{-2}.
\end{eqnarray} 

In the setting of Spin Geometry, a problem similar to the Yamabe problem has been studied in several works of Ammann (see \cite{habilbernd}, \cite{amm1}), and Ammann, Humbert and others (see \cite{amm3}, \cite{amm4}). The starting point of all these works is the Hijazi inequality (\cite{hijazi:86}, \cite{hijazi:91}) which links the first eigenvalue of two elliptic differential operators: the conformal Laplacian $L_g$ and the Dirac operator $D_g$. Hijazi's result can be stated as follow:
\begin{eqnarray}\label{hijazi}
\lambda^2_1(g)\,Vol(M,g)^{\frac{2}{n}}\geq \frac{n}{4(n-1)}Y(M,[g]),
\end{eqnarray}

where $\lambda_1(g)$ denotes the first eigenvalue of the Dirac operator $D_g$. Thereafter, Ammann studies a spin conformal invariant defined by:
\begin{eqnarray}\label{sci}
\bhl:=\underset{\ov{g}\in[g]}{\inf}\lambda_1(\ov{g})Vol(M,\ov{g})^{\frac{1}{n}}
\end{eqnarray}

and points out that studying critical metrics for this invariant involves similar analytic problems to those appearing in the Yamabe problem. Indeed, finding a critical metric of (\ref{sci}) is equivalent to prove the existence of a smooth spinor field $\varphi$ minimizing the functional defined by:
\begin{eqnarray}\label{ahm}
\mathcal{F}_g(\psi)=\frac{\Big(\int_M|D_g\psi|^{\frac{2n}{n+1}}dv(g)\Big)^{\frac{n+1}{n}}}{\big|\int_M\<D_g\psi,\psi\>dv(g)\big|},
\end{eqnarray}

with the corresponding Euler-Lagrange equation given by:
\begin{eqnarray}\label{diraceq}
D_g\varphi=\lambda_{\min}(M,[g],\sigma)|\varphi|^{\frac{2}{n-1}}\varphi.
\end{eqnarray}

In \cite{amm1}, the author observes that a standard variational approach does not yield to the existence of such minimizers. Indeed, the Sobolev inclusion involved in this method is precisely the one for which the compacity is lost in the Reillich-Kondrakov theorem. The argument to overcome this problem is similar to the one used in the Yamabe problem. In fact, one can prove the existence of a smooth solution of Equation~(\ref{diraceq}), but this solution can be trivial (that is identically zero). So one might be able to find a criterion which prevents this situation. It is now important to note that an inequality similar to (\ref{aubinin}) holds in the spinorial setting (see \cite{ammann} and \cite{amm3}), namely:
\begin{eqnarray}\label{large}
\lambda_{\min}(M,[g],\sigma)\leq\lambda_{\min}(\mathbb{S}^n,[\gst],\sigma_{\mathrm{st}})=\frac{n}{2}\omega_n^{\frac{1}{n}}={\sqrt{\frac{n}{n-2}}}K(n,2)^{-1},
\end{eqnarray}

where $(\mathbb{S}^n,\gst,\sigma_{\mathrm{st}})$ stands for the $n$-dimensional sphere equipped with its standard Riemannian metric $\gst$ and its standard spin structure $\sigma_{\mathrm{st}}$. The criterion obtained by Ammann in \cite{amm1} is tightly related to the one involved in the Yamabe problem since he shows that if inequality $(\ref{large})$ is strict then the spinor field solution of (\ref{diraceq}) is non trivial (compare with (\ref{ys})). 

In this paper, we study a more general nonlinear equation involving the Dirac operator (since it also includes Ammann's result in the case of invertible Dirac operator). This equation is closely related to the problem of conformal immersion of a manifold as a hypersurface in a manifold carrying a parallel spinor (see \cite{manuamed} for example). The proof we give here lies on a Sobolev-type inequality for the Dirac operator. It emphasizes in particular that the same kind of questions of those arising from the Yamabe problem can be studied in the context of Spin Geometry.\\

{\bf Acknowledgements:} I would like to thank Oussama Hijazi and Emmanuel Humbert for their encouragements, support and remarks on previous versions of this paper. I am also very grateful to Bernd Ammann for his remarks and suggestions. Finally, I would thank the Mathematical Institute of Neuch\^atel for his financial support.

%%%%%%%%%%%%%%%%%%%%%%%%%%%%%%%%%%%%%%%%%%%%%%%%%%%%%%%%%%%%%%%%%%%%%%%%%%%%%%%%%%%%%%%%%%%%%%%%%%%%%%%%

\section{Geometric and Analytic preliminaries}

%%%%%%%%%%%%%%%%%%%%%%%%%%%%%%%%%%%%%%%%%%%%%%%%%%%%%%%%%%%%%%%%%%%%%%%%%%%%%%%%%%%%%%%%%%%%%%%%%%%%%%%%

%%%%%%%%%%%%%%%%%%%%%%%%%%%%%%%%%%%%%%%%%%%%%%%%%%%%%%%%%%%%%%%%%%%%%%%%%%%%%%%%%%%%%%%%%%%%%%%%%%%%%%%%%

\subsection{Geometric preliminaries}

%%%%%%%%%%%%%%%%%%%%%%%%%%%%%%%%%%%%%%%%%%%%%%%%%%%%%%%%%%%%%%%%%%%%%%%%%%%%%%%%%%%%%%%%%%%%%%%%%%%%%%%%%

In this paragraph, we recall briefly some basic facts on Spin Geometry. For more details, we refer to \cite{fried} or \cite{lawson.michelson:89} for example. Let $(M^n,g,\sigma)$ be an $n$-dimensional compact Riemannian manifold equipped with a spin structure denoted by $\sigma$. It is well-known that on such a manifold one can construct a complex vector bundle of rank $2^{[\frac{n}{2}]}$ denoted by $\Sigma_g(M)$, called the complex spinor bundle. This bundle is naturally endowed with the spinorial Levi-Civita connection $\nabla$, a pointwise Hermitian scalar product $\<.,.\>$ and a Clifford multiplication $``."$. There is also a natural elliptic differential operator of order one acting on sections of this bundle, the Dirac operator. This operator is locally given by:
\begin{eqnarray*}
D_g\varphi=\sum_{i=1}^n e_i\cdot\na_{e_i}\varphi,
\end{eqnarray*}

for all $\varphi\in\Ga\big(\Sigma_g(M)\big)$ and where $\{e_1,...,e_n\}$ is a local $g$-orthonormal frame of the tangent bundle. It defines a self-adjoint operator whose spectrum is constituted of an unbounded sequence of real numbers. Estimates on the spectrum of the Dirac operator has been and is again the main subject of several works (a non exhaustive list is \cite{friedrich}, \cite{hijazi:86} or \cite{baer:92}). As pointed out in the introduction, a key result for the following of this paper is the Hijazi inequality. More precisely, Hijazi gives an inequality which links the squared of the first eigenvalue of the Dirac operator with the first eigenvalue of the conformal Laplacian. The proof of this inequality relies on the famous Schr\"odinger-Lichnerowicz formula (see \cite{cours} for example) and on the conformal covariance of the Dirac operator. In fact, if $\overline{g}\in [g]$, there is a canonical identification between the spinor bundle over $(M,g)$ with the one over $(M,\overline{g})$ (see \cite{hitchin:74} or \cite{hijazi:86}). This identification will be denoted by:
\begin{equation}
\begin{array}{ccc}\label{cc}
\Sigma_{g}(M) & \longrightarrow & \Sigma_{\overline{g}}(M) \\
\varphi & \longmapsto & \overline{\varphi}.
\end{array}
\end{equation}

Under this isomorphism, one can relate the Dirac operators $D_g$ and $D_{\ov{g}}$ acting respectively on $\Sigma_{g}(M)$ and $\Sigma_{\overline{g}} (M)$. Indeed, if $\ov{g}=e^{2u}g$ where $u$ is a smooth function, then:
\begin{eqnarray}\label{diraconforme}
D_{\ov{g}}\,\ov{\varphi}=e^{-\frac{n+1}{2}u}\overline{D_g(e^{\frac{n-1}{2}u}\varphi)},
\end{eqnarray} 

for all $\varphi\in\Gamma\big(\Sigma_{g} (M)\big)$.

%%%%%%%%%%%%%%%%%%%%%%%%%%%%%%%%%%%%%%%%%%%%%%%%%%%%%%%%%%%%%%%%%%%%%%%%%%%%%%%%%%%%%%%

\subsection{Analytic preliminaries}

%%%%%%%%%%%%%%%%%%%%%%%%%%%%%%%%%%%%%%%%%%%%%%%%%%%%%%%%%%%%%%%%%%%%%%%%%%%%%%%%%%%%%%%

In this section we give some well-known facts on Sobolev spaces on spinors and on the analysis of differential equations involving the Dirac operator. In the following, we assume that $(M^n,g)$ is an $n$-dimensional compact Riemannian spin manifold ($n\geq 2$) such that the Dirac operator is invertible.\\

We let $L^q:=L^q\big(\Sigma_g(M)\big)$, the space of spinors $\psi\in\Gamma\big(\Sigma_g(M)\big)$ such that:
\begin{eqnarray*}
||\psi||_q:=\Big(\int_M|\psi|^qdv(g)\Big)^{\frac{1}{q}}
\end{eqnarray*}

is finite. The Sobolev space $H_1^q:=H_1^q\big(\Sigma_g(M)\big)$ is defined as being the completion of the space of smooth spinor fields with respect to the norm:
\begin{eqnarray}
||\varphi||_{1,q}:=||\nabla\varphi||_q+||\varphi||_q.
\end{eqnarray}

However, since our problem involves the Dirac operator, it would be more convenient if we could consider the following equivalent norm:
\begin{lemma}\label{newnorm1}
The map:
\begin{eqnarray}\label{newnorm}
\psi\mapsto||hD_g\psi||_q
\end{eqnarray}

defines a norm equivalent to the $H_1^q$-norm for every smooth positive function $h$ on $M$.
\end{lemma}

{\it Proof:}
From the definition of (\ref{newnorm}) it is clear that this map defines a norm on the space of smooth spinors which is equivalent to the norm defined by:
\begin{eqnarray*}
\psi\mapsto||D_g\psi||_q.
\end{eqnarray*} 

Now we show that this norm is equivalent to the $H_1^q$-norm. In fact, for any smooth spinor field $\psi$, the Cauchy-Schwarz inequality yields:
\begin{eqnarray*}
|D_g\psi|^2\leq n|\nabla\psi|^2
\end{eqnarray*}

which implies the existence of a positive constant $C_1>0$ such that:
\begin{eqnarray*}
||D_g\psi||_q\leq C_1\big(||\nabla\psi||_q+||\psi||_q\big).
\end{eqnarray*} 

On the other hand, with the help of pseudo-differential operators (see the proof of Lemma \ref{pdo}), it is not difficult to see that there also exists another positive constant $C_2>0$ such that:
\begin{eqnarray*}
\big(||\nabla\psi||_q+||\psi||_q\big)\leq C_2||D_g\psi||_q,
\end{eqnarray*} 

which concludes the proof of this lemma.
\hfill $\square$\\
 
Using this result and the fact that the Sobolev space $H^q_1$ is defined as the completion of the space of smooth spinors with respect to the $H_1^q$-norm, it is clear that one can consider the Sobolev space as defined independently from one of the three preceding norms. It will provide a very useful tool to solve the nonlinear equation studied in this paper. A natural way to prove the existence of solutions for this kind of equation is the variational approach. It consists of minimizing a certain functional defined on an adapted Sobolev space and then to apply the machinery of Sobolev-Kondrakov embedding theorems, Schauder estimates and a-priori elliptic estimates. Here we will use this method, and we refer to the works of Ammann (\cite{habilbernd} and \cite{amm1}) for proofs of all these results in the setting of Spin Geometry. However, for clarity, we prove the following result which will be of great help in the next section:  
\begin{lemma}\label{pdo}
If the Dirac operator is invertible then there exists a constant $C>0$ such that for all $\varphi\in H^{q}_1$ we have:
\begin{eqnarray*}
||\varphi||_{p}\leq C ||D_g\varphi||_{q},
\end{eqnarray*}

where $p^{-1}+q^{-1}=1$ and $2\leq p<\infty$.
\end{lemma}

{\it Proof:} We show that the operator:
\begin{eqnarray*}
D_g^{-1}:L^{q}\longrightarrow L^{p}
\end{eqnarray*}

defines a continuous map. Since $D_g^{-1}$ is a pseudo-differential operator of order $-1$ the operator $(Id+\nabla^*\nabla)^{\frac{1}{2}}D_g^{-1}$ is a pseudo-differential operator of order zero hence (see \cite{tayl}) a bounded operator from $L^s$ to $L^s$ for all $s>1$. Thus if $\varphi\in L^q$: 
\begin{eqnarray*}
(Id+\nabla^*\nabla)^{\frac{1}{2}}D_g^{-1}\varphi\in L^q,
\end{eqnarray*}

and the spinor field $D_g^{-1}\varphi$ is in the Sobolev space $H^{q}_1$ which is continuously embedded in $L^{p}$ (using the Sobolev embedding theorem). Then there exists a positive constant $C>0$ such that:
\begin{eqnarray*}
||D_g^{-1}\varphi||_{p}\leq C ||\varphi||_{q},
\end{eqnarray*}

and this concludes the proof.
\hfill$\square$\\

\begin{remark}\label{re1}
\begin{enumerate}
\item For $q=q_D=2n/(n+1)$ and $p_D$ such that $p_D^{-1}+q_D^{-1}=1$ the quotient:
\begin{eqnarray*}
\mathcal{C}_g(\varphi)=\frac{||D_g\varphi||_{q_D}}{||\varphi||_{p_D}}
\end{eqnarray*}

\noindent is invariant under a conformal change of metric, that is:
\begin{eqnarray}\label{covconfepd}
\mathcal{C}_{\ov{g}}(h^{-\frac{n-1}{2}}\varphi)=\mathcal{C}_g(\varphi)
\end{eqnarray}
 
\noindent for all $\varphi\in H^{q_D}_1$ and for $\ov{g}=h^2g\in [g]$. Indeed, an easy computation using the canonical identification (\ref{cc}) between $\Sigma_gM$ and $\Sigma_{\ov{g}}M$ and the formula~(\ref{diraconforme}) which relates $D_g$ and $D_{\ov{g}}$, leads to ${\rm{(\ref{covconfepd})}}$.\\

\item On the $n$-dimensional sphere $(\mathbb{S}^n,g_{{\rm{st}}},\sigma_{{\rm st}})$ endowed with its standard spin structure $\sigma_{{\rm st}}$, the Dirac operator is invertible since the scalar curvature is positive. Then using Lemma~\ref{pdo}, there exists a constant $C>0$ such that for all $\Phi\in H^{q_D}_1$:
\begin{eqnarray*}
||\Phi||_{p_D}\leq C ||D^{\mathbb{S}^n}\Phi||_{q_D}.
\end{eqnarray*}

\noindent Moreover, since the standard sphere $(\mathbb{S}^n\setminus\{q\},\gst)$ (where $q\in\mathbb{S}^n$) is conformally isometric to the Euclidean space $(\mathbb{R}^n,\xi)$, we conclude that for all $\psi\in\Gamma_c\big(\Sigma_{\xi}(\mathbb{R}^n)\big)$:
 \begin{eqnarray*}
||\psi||_{p_D}\leq C ||D_\xi\psi||_{q_D}
\end{eqnarray*}

\noindent where $\Gamma_c\big(\Sigma_{\xi}(\mathbb{R}^n)\big)$ denotes the space of smooth spinor fields over $(\mathbb{R}^n,\xi)$ with compact support.
\end{enumerate}
\end{remark}

%%%%%%%%%%%%%%%%%%%%%%%%%%%%%%%%%%%%%%%%%%%%%%%%%%%%%%%%%%%%%%%%%%%%%%%%%%%%%%%%%%%%%%%%%%%%%%%%%%%%%%%%

\section{The Sobolev Inequality}\label{sectionsobolev}

%%%%%%%%%%%%%%%%%%%%%%%%%%%%%%%%%%%%%%%%%%%%%%%%%%%%%%%%%%%%%%%%%%%%%%%%%%%%%%%%%%%%%%%%%%%%%%%%%%%%%%%%

In this section, we prove a Sobolev inequality in the spinorial setting. The classical Sobolev inequality $S(A,B)$ shows in particular that the Sobolev space of functions $H_1^2$ is continuously embedded in $L^{\frac{2n}{n-2}}$. Here one could interpret our result as the inequality involved in the continuous embedding:
\begin{eqnarray*}
H_1^{2n/(n+1)}\hookrightarrow H_{1/2}^2
\end{eqnarray*}

where $H_{1/2}^2$ is defined as the completion of the space of smooth spinors with respect to the norm:
\begin{eqnarray*}
||\psi||_{\frac{1}{2},2}:=\sum_{i}|\lambda_{i}|^{\frac{1}{2}}|A_i|^2.
\end{eqnarray*}

Here $\psi=\sum_iA_i\psi_i$ is the decomposition of any smooth spinor in the spectral resolution $\{\lambda_i;\psi_i\}$ of $D_g$ (see \cite{habilbernd}).\\

Let's first examine the case of the sphere which is the starting point of the inequality we want to prove. In fact, it is quite easy to compute that the invariant defined by (\ref{sci}) on the sphere is: 
\begin{eqnarray}\label{AHM}
\lambda_{\min}(\mathbb{S}^n,[\gst],\sigma_{\mathrm{st}})=\underset{\psi\neq 0}{\inf}\, \frac{\big(\int_{\mathbb{S}^n}|D^{\mathbb{S}^n}\psi|^{\frac{2n}{n+1}}dv(\gst)\big)^{\frac{n+1}{n}}}{\big|\int_{\mathbb{S}^n}\<D^{\mathbb{S}^n}\psi,\psi\>dv(\gst)\big|}=\frac{n}{2}\omega_n^{\frac{1}{n}}. 
\end{eqnarray}

The proof of this fact relies on the Hijazi inequality (\ref{hijazi}) and on the existence of real Killing spinors on the round sphere (see \cite{sgutt}). Thus using the conformal covariance of (\ref{AHM}) and the fact that the sphere (minus a point) is conformally isometric to the Euclidean space, we can conclude that:
\begin{eqnarray}\label{inre}
\Big|\int_{\mathbb{R}^n}\<D_\xi\psi,\psi\>dx\Big|\leq \lambda_{\min}(\mathbb{S}^n,[\gst],\sigma_{\mathrm{st}})^{-1}\Big(\int_{\mathbb{R}^n}|D_\xi\psi|^{\frac{2n}{n+1}}dx\Big)^{\frac{n+1}{n}}
\end{eqnarray}

for all $\psi\in\Gamma_c\big(\Sigma_{\xi}(\mathbb{R}^n)\big)$. With this in mind, we can now state the main result of this section:
\begin{theorem}\label{sli}
Let $(M^n,g,\sigma)$ be an $n$-dimensional closed compact Riemannian spin manifold and suppose that the Dirac operator is invertible. Then for all $\varepsilon>0$, there exists a constant $B_{\varepsilon}$ such that:
\begin{eqnarray}\label{sli1}
\Big|\int_{M}\<D_g\varphi,\varphi\>dv(g)\Big|\leq \Big(K(n)+\varepsilon\Big)\Big(\int_{M}|D_g\varphi|^{\frac{2n}{n+1}}dv(g)\Big)^{\frac{n+1}{n}}+B_{\varepsilon}\Big(\int_{M}|\varphi|^{\frac{2n}{n+1}}dv(g)\Big)^{\frac{n+1}{n}}
\end{eqnarray}

for all $\varphi\in H^{\frac{2n}{n+1}}_1$ and where 
\begin{eqnarray*}
K(n):=\lambda_{\min}(\mathbb{S}^n,[\gst],\sigma_{\mathrm{st}})^{-1}=\sqrt{\frac{n-2}{n}}K(n,2)=\frac{2}{n}\omega_{n}^{-\frac{1}{n}}.
\end{eqnarray*}
\end{theorem}

\noindent In order to prove (\ref{sli1}), we need some well-known technical results which are summarized in the following lemma:
\begin{lemma}\label{technique}
Let $(a_i)_{1\leq i\leq N_0}\subset\mathbb{R}^+$ $(N_0\in\mathbb{N}^*)$, $p\in [0,1]$ and $q\geq1$. The following identities hold:
\begin{enumerate}

\item $(\sum_{i=1}^{N_0} a_i)^p\leq \sum_{i=1}^{N_0} a_i^p$\\

\item $\sum_{i=1}^{N_0} a_i^q\leq(\sum_{i=1}^{N_0} a_i)^q$\\

\item $\forall\varepsilon>0,\,\exists\,C_{\varepsilon}>0,\,\forall a,b\geq0 \;:\;(a+b)^p\leq (1+\varepsilon)a^p+C_{\varepsilon}b^p$\\

\item For all functions $f_1,\cdots,f_r:M\rightarrow [0,\infty[$, we have:
\begin{eqnarray*}
\sum_{i=1}^r\Big(\int_{M}f^p_i dv(g)\Big)^{\frac{1}{p}}\leq\Big(\int_{M}\big(\sum_{i=1}^r f_i\big)^p dv(g)\Big)^{\frac{1}{p}}.
\end{eqnarray*}
\end{enumerate}
\end{lemma}

We can now give the proof of Inequality (\ref{sli1}).\\

{\it Proof of Theorem~\ref{sli}:} Let $x\in M$ and $\varepsilon>0$. Let $U$ (resp. $\V$) be a neighbourhood of $x\in M$ (resp. $0\in\mathbb{R}^n$) such that the exponential map
\begin{eqnarray*}
exp_x : V\subset\mathbb{R}^n\longrightarrow U\subset M
\end{eqnarray*}

is a diffeomorphism. Then we can identify the spinor bundle over $(U,g)$ with the one over $(V,\xi)$ that is there exists a map:
\begin{eqnarray}\label{id}
\tau : \Sigma_g(U)\longrightarrow\Sigma_\xi(V)
\end{eqnarray}

which is a fiberwise isometry (see \cite{bourguignon}). Moreover the Dirac operators $D_g$ and $D_\xi$ (acting respectively on $\Sigma_g(U)$ and $\Sigma_\xi(V)$) are related by the formula:
\begin{eqnarray}\label{idd}
D_g\varphi(y)=\tau^{-1}\Big(D_\xi\big(\tau(\varphi)\big)\big(exp^{-1}_x(y)\big)\Big)+\rho(\varphi)(y)
\end{eqnarray}

for all $y\in U$ and where $\rho(\varphi)\in\Gamma\big(\Sigma_g(U)\big)$ is a smooth spinor such that $|\rho(\varphi)|\leq \varepsilon|\varphi|$. Now since $M$ is compact, we choose a finite sequence $(x_i)_{1\leq i\leq N_0}\subset M$ and a finite cover $(U_i)_{1\leq i\leq N_0}$ of $M$ (where $U_i$ is a neighbourhood of $x_i\in M$) such that there exists open sets $(V_i)_{1\leq i\leq N_0}$ of $0\in\mathbb{R}^n$ and applications $\tau_i$ such that (\ref{id}) and (\ref{idd}) are fulfilled. Moreover without loss of generality, we can assume that:
\begin{eqnarray*}
\frac{1}{1+\varepsilon}\xi\leq g\leq (1+\varepsilon)\xi
\end{eqnarray*}

as symmetric bilinear forms and consequently the volume forms satisfy:
\begin{eqnarray}\label{volume}
\frac{1}{(1+\varepsilon)^{n/2}}dx\leq dv(g)\leq (1+\varepsilon)^{n/2}dx.
\end{eqnarray}

Let $(\eta_i)_{1\leq i\leq N_0}$ be a smooth partition of unity subordinate to the covering $(U_i)_{1\leq i\leq N_0}$, in other words $\eta_i$ satisfies:
$$\left\lbrace
\begin{array}{l}
{\rm supp\,}(\eta_i)\subset U_i\\
0\leq\eta_i\leq 1\\
\sum_{i=1}^{N_0}\eta_i=1.
\end{array}\right.$$\\

For $\varphi\in\Gamma\big(\Sigma_g(M)\big)$, we write:
\begin{eqnarray*}
(LHS):=\Big|\int_{M}\<D_g\varphi,\varphi\>dv(g)\Big| & = & \Big|\sum_{i=1}^{N_0}\int_{M}\<\sqrt{\eta_i}D_g(\varphi),\sqrt{\eta_i}\varphi\>dv(g)\Big|\\
& =& \Big|\sum_{i=1}^{N_0}\int_{M}\<D_g(\sqrt{\eta_i}\varphi),\sqrt{\eta_i}\varphi\>dv(g)\Big|
\end{eqnarray*}

since $(LHS)$ is real and ${\rm{Re}}\<d(\sqrt{\eta_i})\cdot\varphi,\varphi\>=0$. Inequality~(\ref{volume}) leads to:
\begin{eqnarray*}
(LHS)\leq\big(1+\varepsilon\big)^{\frac{n}{2}}\sum_{i=1}^{N_0}\Big|\int_{\mathbb{R}^n}\<\tau_i\big(D_g(\sqrt{\eta_i}\varphi)\big),\tau_i\big(\sqrt{\eta_i}\varphi\big)\>dx\Big|
\end{eqnarray*}

and using formula~(\ref{idd}), we can write:
\begin{eqnarray*}
(LHS)\leq\big(1+\varepsilon\big)^{\frac{n}{2}}\sum_{i=1}^{N_0}\Big(\Big|\int_{\mathbb{R}^n}\<D_\xi\big(\tau_i(\sqrt{\eta_i}\varphi)\big),\tau_i(\sqrt{\eta_i}\varphi)\>dx\Big|+C\int_{\mathbb{R}^n}|\tau_i(\sqrt{\eta_i}\varphi)|^2dx\Big).
\end{eqnarray*}

On the other hand, since $\tau_i(\sqrt{\eta_i}\varphi)\in\Gamma_c\big(\Sigma_{\xi}(\mathbb{R}^n)\big)$, Inequality~(\ref{inre}) gives:
\begin{eqnarray*}
(LHS)\leq\big(1+\varepsilon\big)^{\frac{n}{2}}\sum_{i=1}^{N_0}\Big(K(n)\big(\int_{\mathbb{R}^n}|D_\xi\big(\tau_i(\sqrt{\eta_i}\varphi)\big)|^{\frac{2n}{n+1}}dx\big)^{\frac{n+1}{n}}+C\int_{\mathbb{R}^n}|\tau_i(\sqrt{\eta_i}\varphi)|^2dx\Big).
\end{eqnarray*}

Now note that with the help of $(4)$ of Lemma~\ref{technique} and since $n/(n+1)\leq 1$, it follows that:
\begin{eqnarray*}
\sum_{i=1}^{N_0}\Big(\int_{\mathbb{R}^n}|D_\xi\big(\tau_i(\sqrt{\eta_i}\varphi)\big)|^{\frac{2n}{n+1}}dx\Big)^{\frac{n+1}{n}} \leq\Big(\int_{\mathbb{R}^n}\big(\sum_{i=1}^{N_0}|D_\xi\big(\tau_i(\sqrt{\eta_i}\varphi)\big)|^2\big)^{\frac{n}{n+1}}dx\Big)^{\frac{n+1}{n}}.
\end{eqnarray*}

Using $(3)$ of Lemma~\ref{technique}, we finally get:
\begin{eqnarray}\label{dirsob}
(LHS)^{\frac{n}{n+1}} & \leq & (1+\varepsilon)^{\frac{(n+1)^2+1}{2(n+1)}} K(n)^{\frac{n}{n+1}} {\bf A}+ (1+\varepsilon)^{\frac{n^2}{2(n+1)}} C_{\varepsilon}{\bf B}
\end{eqnarray}

where:
\begin{eqnarray*}
{\bf A} & = & \int_{\mathbb{R}^n}\Big(\sum_{i=1}^{N_0}|D_\xi\big(\tau_i(\sqrt{\eta_i}\varphi)\big)|^2\Big)^{\frac{n}{n+1}}dx, \textrm{ and }\\
{\bf B} & = & \Big(\sum_{i=1}^{N_0}\int_{\mathbb{R}^n}|\tau_i(\sqrt{\eta_i}\varphi)|^2dx\Big)^{\frac{n}{n+1}}.
\end{eqnarray*}

We now give an estimate of ${\bf A}$. If we let $\gamma_{i}(\varphi)=d(\sqrt{\eta_i})\cdot\varphi-\rho_{i}(\varphi)$, then:
\begin{eqnarray*}
\sum_{i=1}^{N_0}  |D_\xi\big(\tau_i(\sqrt{\eta_i}\varphi)\big)|^2   & = & \sum_{i=1}^{N_0}  |D_g(\sqrt{\eta_i}\varphi)-\rho_i(\varphi)|^2\\
& = & \sum_{i=1}^{N_0} |\sqrt{\eta_i}D_g\varphi+\gamma_{i}(\varphi)|^2
\end{eqnarray*}

and the Minkowski's inequality leads to:
\begin{eqnarray*}
\sum_{i=1}^{N_0} |\sqrt{\eta_i}D_g\varphi+\gamma_{i}(\varphi)|^2  & \leq & \Big(\big(\sum_{i=1}^{N_0}|\sqrt{\eta_i}D_g\varphi|^2\big)^{\frac{1}{2}}+\big(\sum_{i=1}^{N_0}|\gamma_{i}(\varphi)|^2\big)^{\frac{1}{2}}\Big)^2\\
& \leq & \big(|D_g\varphi|+C|\varphi|\big)^2 \qquad(\text{using } (1)\text{ of Lemma~\ref{technique}}).
\end{eqnarray*}

Thus we have shown that:
\begin{eqnarray*}
{\bf A} & \leq & (1+\varepsilon)^{\frac{n}{2}}\int_{M}\big(|D_g\varphi|^2+C|\varphi|^2+C|D_g\varphi|\,|\varphi|\big)^{\frac{n}{n+1}}dv(g),
\end{eqnarray*}

and with $(2)$ of Lemma~\ref{technique}, we get:
\begin{eqnarray*}
{\bf A} & \leq & (1+\varepsilon)^{\frac{n}{2}}\Big(\int_{M}\big(|D_g\varphi|^{\frac{2n}{n+1}}dv(g)+C\int_{M}|\varphi|^{\frac{2n}{n+1}}dv(g)+C\int_{M}|D_g\varphi|^{\frac{n}{n+1}}|\varphi|^{\frac{n}{n+1}}dv(g)\Big).
\end{eqnarray*}

Then we apply the Cauchy-Schwarz inequality in the last term of the preceding inequality:
\begin{eqnarray*}
\int_{M}|D_g\varphi|^{\frac{n}{n+1}}|\varphi|^{\frac{n}{n+1}}dv(g) & \leq & \Big(\int_{M}|D_g\varphi|^{\frac{2n}{n+1}}dv(g)\Big)^{\frac{1}{2}}\Big(\int_{M}|\varphi|^{\frac{2n}{n+1}}dv(g)\Big)^{\frac{1}{2}}
\end{eqnarray*}

and next we use the Young inequality:
\begin{eqnarray*}
\int_{M}|D_g\varphi|^{\frac{n}{n+1}}|\varphi|^{\frac{n}{n+1}}dv(g) & \leq &\frac{\varepsilon^2}{2}\int_{M}|D_g\varphi|^{\frac{2n}{n+1}}dv(g)+\frac{1}{2\varepsilon^2}\int_{M}|\varphi|^{\frac{2n}{n+1}}dv(g).
\end{eqnarray*}

Finally, we have:
\begin{eqnarray*}
{\bf A} \leq 
(1+\varepsilon)^{\frac{n}{2}}\Big((1+\frac{\varepsilon^2}{2})\int_{M}|D_g\varphi|^{\frac{2n}{n+1}}dv(g)+C_\varepsilon\int_{M}|\varphi|^{\frac{2n}{n+1}}dv(g)\Big).
\end{eqnarray*}

Now we estimate ${\bf B}$ in Inequality~(\ref{dirsob}). H\"older's inequality gives:
\begin{eqnarray*}
\int_{\mathbb{R}^n}|\tau_i(\sqrt{\eta_i}\varphi)|^2dx\leq \Big(\int_{\mathbb{R}^n}|\tau_i(\sqrt{\eta_i}\varphi)|^{\frac{2n}{n-1}}dx\Big)^{\frac{n-1}{2n}}\Big(\int_{\mathbb{R}^n}|\tau_i(\sqrt{\eta_i}\varphi)|^{\frac{2n}{n+1}}dx\Big)^{\frac{n+1}{2n}}
\end{eqnarray*}

and using $(1)$ of Lemma~\ref{technique} and the preceding inequality lead to:
\begin{eqnarray*}
{\bf B} & \leq & \sum_{i=1}^{N_0} 
\Big(\int_{\mathbb{R}^n}|\tau_i(\sqrt{\eta_i}\varphi)|^{\frac{2n}{n-1}}dx\Big)^{\frac{n-1}{2(n+1)}}\Big(\int_{\mathbb{R}^n}|\tau_i(\sqrt{\eta_i}\varphi)|^{\frac{2n}{n+1}}dx\Big)^{\frac{1}{2}}.
\end{eqnarray*}

With the help of $(2)$ of Remark~\ref{re1}, there exists a constant $C>0$ such that:
\begin{eqnarray*}
{\bf B} & \leq & C \sum_{i=1}^{N_0} 
\Big(\int_{\mathbb{R}^n}|D_\xi\big(\tau_i(\sqrt{\eta_i}\varphi)\big)|^{\frac{2n}{n+1}}dx\Big)^{\frac{1}{2}}\Big(\int_{\mathbb{R}^n}|\tau_i(\sqrt{\eta_i}\varphi)|^{\frac{2n}{n+1}}dx\Big)^{\frac{1}{2}}.
\end{eqnarray*}

On the other hand, the Young inequality gives:
\begin{eqnarray*}
{\bf B} & \leq & C\varepsilon^2\sum_{i=1}^{N_0} \int_{\mathbb{R}^n}|D_\xi\big(\tau_i(\sqrt{\eta_i}\varphi)\big)|^{\frac{2n}{n+1}}dx+\frac{C}{\varepsilon^2}\int_{\mathbb{R}^n}|\tau_i(\sqrt{\eta_i}\varphi)|^{\frac{2n}{n+1}}dx.
\end{eqnarray*}

and it is easy to see that:
\begin{eqnarray*}
\frac{C}{\varepsilon^2}\sum_{\i=1}^{N_0}\int_{\mathbb{R}^n}|\tau_i(\sqrt{\eta_i}\varphi)|^{\frac{2n}{n+1}}dx & \leq & C_\varepsilon (1+\varepsilon)^{\frac{n}{2}}\int_{M}|\varphi|^{\frac{2n}{n+1}}dv(g).
\end{eqnarray*}

To conclude, an argument similar to the one used in the estimate of ${\bf A}$ shows that:
\begin{eqnarray*}
\sum_{i=1}^{N_0}\int_{\mathbb{R}^n}|D_\xi\big(\tau_i(\sqrt{\eta_i}\varphi)\big)|^{\frac{2n}{n+1}}dx\leq (1+\varepsilon)^{\frac{n}{2}}\int_{M}\big(|D_g\varphi|^{\frac{2n}{n+1}}+C|\varphi|^{\frac{2n}{n+1}}+|D_g\varphi|^{\frac{n}{n+1}}|\varphi|^{\frac{n}{n+1}}\big)dv(g)
\end{eqnarray*}

and the Cauchy-Schwarz inequality and the Young inequality lead to:
\begin{eqnarray*}
{\bf B} & \leq & C\varepsilon^2(1+\varepsilon)^{\frac{n}{2}}\int_{M}|D_g\varphi|^{\frac{2n}{n+1}}dv(g)+C_{\varepsilon}\int_{M}|\varphi|^{\frac{2n}{n+1}}dv(g).
\end{eqnarray*}

Combining the estimates of ${\bf A}$ and ${\bf B}$ in (\ref{dirsob}) gives Inequality~(\ref{sli1}).
\hfill$\square$\\

%%%%%%%%%%%%%%%%%%%%%%%%%%%%%%%%%%%%%%%%%%%%%%%%%%%%%%%%%%%%%%%%%%%%%%%%%%%%%%%%%%%%%%%%%%%%%%%%%%%%%%%%

\section{A nonlinear equation for the Dirac operator}

%%%%%%%%%%%%%%%%%%%%%%%%%%%%%%%%%%%%%%%%%%%%%%%%%%%%%%%%%%%%%%%%%%%%%%%%%%%%%%%%%%%%%%%%%%%%%%%%%%%%%%%%

%%%%%%%%%%%%%%%%%%%%%%%%%%%%%%%%%%%%%%%%%%%%%%%%%%%%%%%%%%%%%%%%%%%%%%%%%%%%%%%%%%%%%%%%%%%%%%%%%%%%%%%%

\subsection{A criterion for the existence of solutions}

%%%%%%%%%%%%%%%%%%%%%%%%%%%%%%%%%%%%%%%%%%%%%%%%%%%%%%%%%%%%%%%%%%%%%%%%%%%%%%%%%%%%%%%%%%%%%%%%%%%%%%%%

As a direct application of Theorem~\ref{sli}, we give a sufficient criterion for the existence of solutions for a nonlinear equation involving the Dirac operator. More precisely, the aim of this section is to prove the following result:
\begin{theorem}\label{existence}
Let $(M^n,g)$ be an $n$-dimensional compact Riemannian spin manifold and let $H$ be a smooth positive function on $M$. If the Dirac operator is invertible and if:
\begin{eqnarray}\label{strictdirac}
\lambda_{\min}<K(n)^{-1}(\max H)^{-\frac{2}{p_D}},
\end{eqnarray}

then there exists a spinor field $\varphi\in C^{1,\alpha}(M)\cap\ C^{\infty}\big(M\setminus \varphi^{-1}(0)\big)$ satisfying the following nonlinear elliptic equation:
\begin{eqnarray}\label{nld}
D_g\varphi= \lambda_{\min} H|\varphi|^{\frac{2}{n-1}}\varphi\quad\text{and}\quad\int_MH|\varphi|^{\frac{2n}{n-1}}dv_g=1.
\end{eqnarray}
\end{theorem}

In the statement of Theorem \ref{existence}, we let for $2\leq q\leq q_D$:
\begin{eqnarray}\label{lambdaq}
\lambda_q=\lambda_q(M,g,\sigma):=\underset{\psi\neq 0}{\inf}\Big\{\frac{\big(\int_{M}H^{-(q/p)}|D_g\psi|^{q}dv(g)\big)^{\frac{2}{q}}}{\big|\int_{M}\<D_g\psi,\psi\>dv(g)\big|}\Big\}=\underset{\psi\neq 0}{\inf}\,\frac{||H^{-(1/p)}D_g\psi||_q^2}{\big|\int_{M}\<D_g\psi,\psi\>dv(g)\big|}
\end{eqnarray}

where the infimum is taken over all $\psi\in H^q_1$ and where $\lambda_{q_D}(M,g,\sigma):=\lambda_{\min}$. In the rest of this section, we will let:
\begin{eqnarray*}
\mathcal{F}_q(\psi)=\mathcal{F}_{g,q}(\psi)=\frac{||H^{-(1/p)}D_g\psi||_q^2}{\big|\int_{M}\<D_g\psi,\psi\>dv(g)\big|}.
\end{eqnarray*}

Here $C^\infty(M)$ (resp. $C^{k,\alpha}(M)$) denotes the space of smooth spinor fields (resp. of spinor fields with finite $(k,\alpha)$-H\"older norm) on $M$ (see \cite{habilbernd}).
 
\begin{remark}\label{re2}
Using Lemma~\ref{pdo}, we have $\lambda_q>0$.
\end{remark}

A standard variational approach to study (\ref{nld}) cannot allow to conclude because of the lack of compactness of the inclusion $H_1^{q_D}$ in $L^{p_D}$. The method we use here consists in proving the existence of solutions for subcritical equations where the compactness of the Sobolev embedding theorem is valid. Then we prove that one can extract a subsequence which converges to a solution of (\ref{nld}). We begin with the existence of solutions for subcritical equations, that is:
\begin{proposition}\label{equationsc}
For all $q\in(q_D,2)$, there exists a spinor field $\varphi_q\in C^{1,\alpha}(M)\cap\ C^{\infty}\big(M\setminus \varphi^{-1}_q(0)\big)$ such that:
\begin{eqnarray*}
D_g\varphi_q=\lambda_{q}H|\varphi_q|^{p-2}\varphi_q\qquad(E_{q})
\end{eqnarray*}

where $p\in\mathbb{R}$ is such that $p^{-1}+q^{-1}=1$. Moreover, we have:
\begin{eqnarray*}
\int_MH|\varphi_q|^{p}dv_g=1.
\end{eqnarray*}
\end{proposition} 

{\it Proof:} The proof of this result is divided into two parts. In a first step, we show that there exists a spinor field $\varphi_q\in H^q_1$ satisfying $(E_{q})$, and then we will show that this solution has the desired regularity. For the rest of this proof, we fix $q\in (q_D,2)$.\\

{\it First step:} We prove the existence of a spinor field $\varphi_q\in  H^q_1$ satisfying $(E_{q})$. First we study the functional defined by:
\begin{eqnarray*}
\mathcal{F}_q:\mathcal{H}_1^q:=\big\{\psi\in H_1^q\,\,/\,\int_{M}\<D_g\psi,\psi\>dv(g)=1\big\}\longrightarrow \mathbb{R}.
\end{eqnarray*}

It is clear that $\mathcal{H}_1^q$ is non empty. Take for example a smooth eigenspinor $\psi_1$ associated to the first positive eigenvalue $\lambda_1>0$ of the Dirac operator and thus $(\lambda_1)^{-(1/2)}||\psi_1||_2^{-1}\psi_1\in\mathcal{H}_1^q$. On the other hand, since $\mathcal{F}_q(\psi)\geq 0$ for all $\psi\in\mathcal{H}^q_1$, we can consider a minimizing sequence $(\psi_i)$ for $\mathcal{F}_q$, that is a sequence such that $\mathcal{F}_q(\psi_i)\rightarrow\lambda_q$ with $(\psi_i)\subset\mathcal{H}_1^{q}$. It is clear that this sequence is bounded in $H^q_1$ and thus there exists a spinor field $\psi_q\in H_1^q$ such that:
\begin{itemize}

\item $\psi_i\rightarrow\psi_q$ strongly in $L^p$ with $p^{-1}+q^{-1}=1$ (by the Reillich-Kondrakov theorem)

\item $\psi_i\rightarrow\psi_q$ weakly in $H_1^q$ (by reflexivity of the Sobolev space $H_1^q$).
 
\end{itemize}

Moreover, we write:
\begin{eqnarray*}
\int_{M}\<D_g\psi_q,\psi_q\>dv(g)=\int_{M}\<D_g\psi_q,\psi_q-\psi_i\>dv(g)+\int_{M}\<D_g\psi_q,\psi_i\>dv(g)
\end{eqnarray*}

and we note that:
\begin{eqnarray*}
\big|\int_{M}\<D_g\psi_q,\psi_q-\psi_i\>dv(g)\big|\leq ||D_g\psi_q||_q ||\psi_q-\psi_i||_p\longrightarrow 0
\end{eqnarray*}

where we used the H\"older inequality and the strong convergence in $L^p$. One can also easily check that the map:
\begin{eqnarray*}
\Phi\longmapsto\int_{M}\<D_g\psi_q,\Phi\>dv(g)
\end{eqnarray*}

defines a continuous linear form on $H^q_1$ and then the weak convergence in $H^q_1$ gives:
\begin{eqnarray*}
\int_{M}\<D_g\psi_q,\psi_q\>dv(g)=1,
\end{eqnarray*}

that is $\psi_q\in\mathcal{H}^q_1$. Once again because of the weak convergence in $H^q_1$ and of Lemma~\ref{newnorm1}, we also have:
\begin{eqnarray*}
||H^{-(1/p)}D_g\psi_q||^2_q\leq\underset{i\rightarrow\infty}{\lim\inf\,}||H^{-(1/p)}D_g\psi_i||^2_q=\lambda_q
\end{eqnarray*}
 
and thus $\lambda_q=\mathcal{F}_q(\psi_q)$. Finally, we proved that there exists $\psi_q\in\mathcal{H}^q_1$ which reaches $\lambda_q$. Now for all smooth spinors $\Phi$, we compute:
\begin{eqnarray*}
\frac{d}{dt}_{|t=0} ||D_g(\psi_q+t\Phi)||_q^2 & = & 2\lambda_q^{\frac{2-q}{2}}\int_{M}\mathrm{Re}\<H^{-(q/p)}|D_g\psi_q|^{q-2}D_g\psi_q,D_g\Phi\>dv(g)
\end{eqnarray*}

and:
\begin{eqnarray*}
\frac{d}{dt}_{|t=0} \int_{M}\mathrm{Re}\<D_g(\psi_q+t\Phi),(\psi_q+t\Phi)\>dv(g) = 2\int_M\mathrm{Re}\<\psi_q,D_g\Phi\>dv(g)
\end{eqnarray*}

which, by the Lagrange multipliers theorem, gives the existence of a real number $\alpha$ such that:
\begin{eqnarray*}
\lambda_q^{\frac{2-q}{2}}\int_{M}\mathrm{Re}\<H^{-(q/p)}|D_g\psi_q|^{q-2}D_g\psi_q,D_g\Phi\>dv(g)= \alpha \int_M\mathrm{Re}\<\psi_q,D_g\Phi\>dv(g).
\end{eqnarray*}

Moreover, since $\psi_q$ is a critical point for $\mathcal{F}_q$, we get $\alpha=\lambda_q$ and thus:
\begin{eqnarray*} 
\int_{M}\<\lambda_q^{q/2}\psi_q-H^{-(q/p)}|D_g\psi_q|^{q-2}D_g\psi_q,D_g\Phi\>dv(g)=0.
\end{eqnarray*}

To sum up, we proved the existence of a spinor field $\psi_q\in\mathcal{H}^q_1$ satisfying weakly the equation:
\begin{eqnarray*}
|D_g\psi_q|^{q-2}D_g\psi_q=\lambda_q^{q/2}H^{q/p}\psi_q.
\end{eqnarray*}

If we let $\phi_q=\lambda_q^{1/2}\psi_q$, we can easily check that $\phi_q\in H^q_1$ satisfies $(E_q)$ (where we used the relations $|\psi_q|=\lambda_q^{-(q/2)}H^{-(q/p)}|D_g\psi_q|^{q/p}$ and $|D_g\psi_q|^{2-q}=(\lambda_q^{q/2}H^{q/p}|\psi_q|)^{p-2}$). On the other hand, since:
\begin{eqnarray*}
\int_{M}\<D_g\psi_q,\psi_q\>dv(g)=1,
\end{eqnarray*}

and since the spinor field $\varphi_q$ is a solution of $(E_q)$, we deduce that:
\begin{eqnarray*}
\int_MH|\varphi_q|^{p}dv_g=1.\\
\end{eqnarray*}

{\it Second step:} We show that $\varphi_q\in C^{1,\alpha}(M)\cap\ C^{\infty}\big(M\setminus \varphi^{-1}_q(0)\big)$. The proof of this result uses the classical ``bootstrap argument''. Indeed, the spinor field $\varphi_q$ is in the Sobolev space $H^q_1$ which is continuously embedded in $L^{p_1}$ with $p_1=nq/(n-q)$, by the Sobolev embedding theorem. The H\"older inequality implies that $H|\varphi_q|^{p-2}\varphi\in L^{p_1/(p-1)}$ and then elliptic a-priori estimates~(see \cite{amm1}) gives $\varphi\in H^{p_1/(p-1)}_1$. Once again, the Sobolev embedding theorem implies that $\varphi_q\in L^{p_2}$ with
\begin{eqnarray*}
p_2=np_1/(n(p-1)-p_1),
\end{eqnarray*}

if $n(p-1)>p_1$ or $\phi_q\in L^s$ for all $s>1$ if $n(p-1)\leq p_1$. Note that since $q>q_D$, we can easily check that $p_2>p_1$ and thus we get a better regularity for the spinor field $\varphi_q$. In fact, if we push further this argument, we can show that $\varphi_q\in L^{p_i}$ for all $i$, where $p_i$ is the sequence of real numbers defined by:
$$p_i:=\left\lbrace
\begin{array}{ll}
\frac{np_{i-1}}{n(p-1)-p_{i-1}}\quad & \text{if}\quad n(p-1)>p_{i-1}\\ \\
+\infty \quad & \text{if}\quad n(p-1)\leq p_{i-1}.
\end{array}
\right.$$

A classical study of this sequence leads to the existence of a rank $i_0\in\mathbb{N}$ such that $p_{i_0}=+\infty$ and thus we can conclude that $\varphi_q\in L^s$ for all $s>1$. The elliptic a-priori estimate gives that $\varphi_q\in H^s_1$ for all $s>1$ and if we apply the Sobolev embedding theorem, one concludes that $\varphi_q\in C^{0,\alpha}(M)$ for $\alpha\in (0,1)$. Hence $f|\varphi_q|^{p-2}\varphi_q\in C^{0,\alpha}(M)$ as well, and the Schauder estimate (see \cite{amm1}) gives $\varphi_q\in C^{1,\alpha}(M)$. It is clear that one can carry on this argument on $M\setminus\varphi_{q}^{-1}(0)$ to obtain $\varphi_q\in C^{\infty}\big(M\setminus\varphi_{q}^{-1}(0)\big)$. 
\hfill$\square$\\

\begin{remark}
If we assume that $p\geq 2$, the regularity of the spinor field $\varphi_q$ can be improved to $C^{2,\alpha}(M)$.
\end{remark}

In the following, we want to prove the existence of a solution of the equation $(E_{q_D})$. However, we cannot argue like in the proof of Proposition~\ref{equationsc} because of the lack of compacity of the embedding $H^{q_D}_1\hookrightarrow L^{p_D}$ which is precisely the one involved in our problem. The idea is to adapt the proof of the Yamabe problem (see for example \cite{lee.parker:87}). Indeed we will prove that one can extract a subsequence from the sequence of solutions $(\varphi_q)$ which converges to a weak solution of Problem~(\ref{nld}) (see Lemma~\ref{partone}). Then in Lemma~\ref{parttwo}, we will get the desired regularity for this solution and finally in Lemma~\ref{partthree}, using Inequality~(\ref{sli1}) of Theorem~\ref{sli}, we will be able to exclude the trivial solution. So we first have:
\begin{lemma}\label{partone}
There exists a sequence $(q_i)$ which tends to $q_D$ and such that the corresponding sequence $(\varphi_{q_i})$, solution of $(E_{q_i})$, converges to a weak solution $\varphi\in H^{q_D}_1$ of (\ref{nld}).
\end{lemma}

{\it Proof:}
It is clear that without loss of generality, we can suppose that the volume of the manifold $(M,g)$ is equal to $1$. Otherwise, because of the conformal covariance of Equation (\ref{nld}), we change the metric with a homothetic one (and so a conformal one). In a similar way, we can also assume (because of a rescaling argument) that the maximum of the function $H$ is equal to $1$. Now we prove that the sequence $(\varphi_q)$ is uniformly bounded in $H^{q_D}_1$. Indeed, since $q\geq q_D$, the H\"older inequality gives:
\begin{eqnarray*}
||H^{-(1/p_D)}D_g\varphi_q||_{q_D}^2\leq ||H^{-(1/p_D)}D_g\varphi_q||_{q}^2.
\end{eqnarray*}

On the other hand, $p\leq p_D$ implies that:
\begin{eqnarray*}
||H^{-(1/p_D)}D_g\varphi_q||_{q_D}^2\leq\lambda^2_q.
\end{eqnarray*}

The variational characterization of $\lambda_q$ and the H\"older inequality directly yield:
\begin{eqnarray*}
||H^{-(1/p_D)}D_g\varphi_q||_{q_D}^2\leq\lambda_q^2\leq\lambda_1^2(\min H)^{-1}
\end{eqnarray*}

and thus we conclude that $(\varphi_q)$ is uniformly bounded in $H^{q_D}_1$. Then there exists a sequence $(q_i)$ which tends to $q_D$ and a spinor field $\varphi\in H^{q_D}_1$ such that:
\begin{itemize}
\item $\varphi_{q_i}\rightarrow\varphi$ weakly in $H_1^{q_D}$ (by reflexivity of the Sobolev space $H_1^{q_D}$)

\item $\varphi_{q_i}\rightarrow\varphi$ a.e. on $M$.
\end{itemize}

Moreover, since $(\varphi_{q_i})$ is bounded in $H^{q_D}_1$, the Sobolev embedding theorem implies that it is bounded in $L^{p_D}$, and so $H|\varphi_{q_i}|^{p_i-2}\varphi_{q_i}$ is bounded in $L^{p_D/(p_i-1)}$. However, since $p_D/(p_D-1)<p_D/(p_i-1)$, the sequence $H|\varphi_{q_i}|^{p_i-2}\varphi_{q_i}$ is also bounded in $L^{p_D/(p_D-1)}$. Using this fact and since 
\begin{eqnarray*}
H|\varphi_{q_i}|^{p_i-2}\varphi_{q_i}\rightarrow H|\varphi|^{p_D-2}\varphi\quad\text{a.e. on } M,
\end{eqnarray*}

we finally get that:
\begin{eqnarray*}
H|\varphi_{q_i}|^{p_i-2}\varphi_{q_i}\rightarrow H|\varphi|^{p_D-2}\varphi\quad\text{weakly in }L^{p_D/(p_D-1)}
\end{eqnarray*}

and so weakly in $L^1$. Now note that for all smooth spinor fields $\Phi$, the map:
\begin{eqnarray*}
\psi\longmapsto\int_{M}\<D_g\psi,D_g\Phi\>dv(g)
\end{eqnarray*}

defines a continuous linear form on $H^{q_D}_1$ and thus by weak convergence in $H^{q_D}_1$, we obtain:
\begin{eqnarray*}
\int_{M}\<D_g\varphi_{q_i},D_g\Phi\>dv(g)\underset{i\rightarrow +\infty}{\longrightarrow}\int_{M}\<D_g\varphi,D_g\Phi\>dv(g).
\end{eqnarray*}

The weak convergence in $L^1$ gives:
\begin{eqnarray*}
\int_{M}\<H|\varphi_{q_i}|^{p_i-2}\varphi_{q_i},D_g\Phi\>dv(g)\underset{i\rightarrow +\infty}{\longrightarrow}\int_{M}\<H|\varphi|^{p_D-2}\varphi,D_g\Phi\>dv(g).
\end{eqnarray*}

Now using the variational characterization~(\ref{lambdaq}) of $\lambda_q$ and the fact that the function:
\begin{eqnarray*}
q\longmapsto||D_g\Phi||_q
\end{eqnarray*}

is continuous, we easily conclude that $q\longmapsto\lambda_q$ is also continuous. Combining all the preceding statements with the fact that $\varphi_{q_i}$ is a solution of $(E_{q_i})$ leads to:
\begin{eqnarray*}
\int_{M}\<D_g\varphi,D_g\Phi\>dv(g)=\lambda_{\min}\int_{M}\<H|\varphi|^{2/(n-1)}\varphi,D_g\Phi\>dv(g),
\end{eqnarray*}

for all smooth spinor fields $\Phi$, that is $\varphi\in H^{q_D}_1$ is a weak solution of (\ref{nld}).
\hfill$\square$\\

We then state a regularity Lemma which is proved in \cite{amm1} and thus we omit the proof here. 
\begin{lemma}\label{parttwo}
The spinor field $\varphi$ given in Lemma~\ref{partone} satisfies $\varphi\in C^{1,\alpha}(M)\cap\ C^{\infty}\big(M\setminus \varphi^{-1}(0)\big)$.
\end{lemma}

As pointed out by Tr\"udinger in the context of the Yamabe problem, one cannot exclude from this step the case where the spinor field $\varphi$, obtained in Lemma~\ref{partone} and \ref{parttwo}, is identically zero. In \cite{amm1}, Ammann proves that if (\ref{strictdirac}) (with $H$ constant) is fulfilled then $\varphi$ is non trivial. We give a similar result for Equation~(\ref{nld}) which generalizes the one of Ammann in the case where the Dirac operator is invertible. The proof we present here is based on the Sobolev-type inequality obtained in Theorem \ref{sli}. More precisely, we get:
\begin{lemma}\label{partthree}
If (\ref{strictdirac}) is satisfied, the spinor $\varphi$ obtained in Lemma~\ref{partone} and \ref{parttwo} is non identically zero and:
\begin{eqnarray*}
\int_MH|\varphi|^{\frac{2n}{n-1}}dv_g=1.
\end{eqnarray*}
\end{lemma}

{\it Proof:}
Let $\varphi_q\in C^{1,\alpha}(M)\cap\ C^{\infty}\big(M\setminus\varphi^{-1}_q(0)\big)$ be a solution of Equation~$(E_q)$, that is:
\begin{eqnarray*}
D_g\varphi_q=\lambda_{q}H|\varphi_q|^{p-2}\varphi_q 
\end{eqnarray*}

and $\int_MH|\varphi_q|^pdv_g=1$ for all $q\in(q_D,2)$ (where $p$ is such that $p^{-1}+q^{-1}=1$). Since $q>q_D$, the H\"older inequality yields:
\begin{eqnarray*}
\Big(\int_M|D_g\varphi_q|^{q_D}dv(g)\Big)^{2/q_D}\leq\big(\max H\big)^{2/p}\Big(\int_M|H^{-(1/p)}D_g\varphi_q|^{q}dv(g)\Big)^{2/q} Vol(M,g)^{2(q-q_D)/(qq_D)}
\end{eqnarray*}

and with the help of $(E_q)$, we get:
\begin{eqnarray*}
\int_M|H^{-(1/p)}D_g\varphi_q|^{q}dv(g)=\lambda^q_q.
\end{eqnarray*}

We finally obtain:
\begin{eqnarray}\label{pasdeposte}
\big(\int_M|D_g\varphi_q|^{q_D}dv(g)\big)^{2/q_D}\leq \big(\max H\big)^{2/p}\lambda^2_q\,Vol(M,g)^{2(q-q_D)/(qq_D)}.
\end{eqnarray}

On the other hand, applying Theorem~\ref{sli} for the spinor fields $\varphi_q$ gives:
\begin{eqnarray*}
\int_{M}\<D_g\varphi_q,\varphi_q\>dv(g)=\lambda_q\leq \Big(K(n)+\varepsilon\Big)\Big(\int_{M}|D_g\varphi_q|^{q_D}dv(g)\Big)^{2/q_D}+B_{\varepsilon}\Big(\int_{M}|\varphi_q|^{q_D}dv(g)\Big)^{2/q_D}
\end{eqnarray*}

where $B_{\varepsilon}>0$ is a positive constant. Using $(\ref{pasdeposte})$ in the preceding inequality leads:
\begin{eqnarray*}
1\leq \Big(K(n)+\varepsilon\Big)\big(\max H\big)^{2/p}\lambda_q\; Vol(M,g)^{2(q-q_D)/(qq_D)}+B_{\varepsilon}\Big(\int_{M}|\varphi_q|^{q_D}dv(g)\Big)^{2/q_D}.
\end{eqnarray*}

Now if $q$ tends to $q_D$, we obtain:
\begin{eqnarray*}
1\leq \Big(K(n)+\varepsilon\Big)\big(\max H\big)^{2/p_{D}}\lambda_{\min}+B_{\varepsilon}\Big(\int_{M}|\varphi|^{q_D}dv(g)\Big)^{2/q_{D}}.
\end{eqnarray*}

However, because of (\ref{strictdirac}), we have:
\begin{eqnarray*}
\big(\max H\big)^{2/p_{D}}\lambda_{\min}\,K(n)<1,
\end{eqnarray*}

which allows to conclude that, for $\varepsilon>0$ small enough, the norm $||\varphi||_{q_D}>0$ and thus $\varphi$ is not identically zero.
\hfill$\square$\\

\begin{remark}
Note that we recover the result of Ammann proved in \cite{amm1} for $H=cste$ (under the assumption that the Dirac operator has a trivial kernel).\\
\end{remark}

%%%%%%%%%%%%%%%%%%%%%%%%%%%%%%%%%%%%%%%%%%%%%%%%%%%%%%%%%%%%%%%%%%%%%%%%%%%%%%%%%%%%%%%%%%%%%%%%%%%%%%%%

\subsection{An upper bound for $\lambda_{\min}$}\label{partwo}

%%%%%%%%%%%%%%%%%%%%%%%%%%%%%%%%%%%%%%%%%%%%%%%%%%%%%%%%%%%%%%%%%%%%%%%%%%%%%%%%%%%%%%%%%%%%%%%%%%%%%%%%

In this section, we prove a general upper bound for $\lambda_{\min}$. Namely, we get:
\begin{theorem}\label{existance}
Let $(M^n,g)$ be an $n$-dimensional compact Riemannian spin manifold with $n\geq 3$. If $H\in C^\infty(M)$ is a smooth positive function on $M$, then the following inequality holds:
\begin{eqnarray*}
\lambda_{\min}\leq K(n)^{-1}(\max_M H)^{-\frac{2}{p_D}}.
\end{eqnarray*}
\end{theorem}

The proof of Theorem~\ref{existance} lies on the construction of an adapted test spinor which will be estimated in the variational characterization of $\lambda_{\min}$. We first note that $\lambda_{\min}$ is invariant under a conformal change of the metric, therefore we can work with any metric within the conformal class of $g$. Indeed, we have:
\begin{proposition}
The number $\lambda_{\min}$ is a conformal invariant of $(M,g)$.
\end{proposition}

{\it Proof:}
We can easily compute that for $\ov{g}=u^2 g\in[g]$ we have:
\begin{eqnarray*}
\mathcal{F}_{\ov{g},q_D}(\ov{\psi})=\mathcal{F}_{g,q_D}(u^{\frac{n-1}{2}}\psi)
\end{eqnarray*}

and then because of the variational characterization (\ref{lambdaq}) of $\lambda_{\min}(M,\ov{g},\sigma)$, its conformal covariance follows directly. 
\hfill$\square$\\

For sake of completeness, we briefly recall the work of Ammann, Grosjean, Humbert and Morel \cite{amm3} which describes in particular the construction of the test-spinor. We first need a trivialization of the spinor bundle given by the Bourguignon-Gauduchon trivialization \cite{bourguignon} which is adapted for our problem. Let $(x_1,...,x_n)$ be the Riemannian normal coordinates given by the exponential map at $p\in M$:
$$\begin{array}{lclc}
\exp_p : & V\subset T_{p}M\simeq\mathbb{R}^n & \longrightarrow & U\subset M \\
 & (x_1,...,x_{n}) & \longmapsto & m.
\end{array}$$

Now if we consider the smooth map $m\mapsto G_m:=\big(g_{ij}(m)\big)$ which associates to any point $m\in U$ the matrix of the coefficients of the metric $g$ at this point in the basis $\{\frac{\pa}{\pa x_1},...,\frac{\pa}{\pa x_n}\}$, then one can find a unique symmetric matrix $B_m:=\big(b^j_i(m)\big)$ (which depends smoothly on $m$) such that $B^2_m=G^{-1}_m$. Thus, at each point $m\in U$ we obtain an isometry between $\mathbb{R}^n$ and the tangent space $T_mM$ defined by:
$$\begin{array}{lcll}
\mathrm{B}_m: & \big( \T_{\exp^{-1}_{p}(m)}V\simeq\mathbb{R}^{n},\xi \big) & \longrightarrow & \big(T_mU,g_m\big) \\
 & (a^1,...,a^n) & \longmapsto & \sum_{i,j}b^j_i(m)a^i\frac{\pa}{\pa x_j}(m).
\end{array}$$

This map induces an identification between the two $SO_n$-principal bundles of orthonormal frames over $(V,\xi)$ and $(U,g)$. Thereafter, this identification can be lifted to the $Spin_n$-principal bundles of spinorial frames over $(V,\xi)$ and $(U,g)$ and then gives an isometry:
$$\begin{array}{clc}
\Sigma_\xi(V) & \longrightarrow & \Sigma_g(U) \\
 \varphi & \longmapsto & \ov{\varphi}.
\end{array}$$
 
This identification has already been used in Section~\ref{sectionsobolev} and was denoted by $\tau$. However, for sake of clarity, we will denote it by $\tau(\varphi):=\ov{\varphi}$ for $\varphi\in\Gamma\big(\Sigma_\xi(V)\big)$. Now let:
\begin{eqnarray*}
e_i:=b^j_i\frac{\partial}{\partial x_j},
\end{eqnarray*}

such that $\{e_1,...,e_n\}$ defines an orthonormal frame of $(TU,g)$. Via the preceding identification, one can relate the Dirac operator acting on $\Sigma_\xi(V)$ with the one acting on $\Sigma_g(U)$. Indeed, if $D_\xi$ and $D_g$ denote those Dirac operators, we have:
\begin{eqnarray}\label{diracrelation}
D_g\overline{\psi}=\overline{D_\xi\psi}+\sum_{i,j=1}^{n}\big(b_i^j-\delta_i^j\big)\ov{\pa_i\cdot\na_{\pa_j}\psi}+{\bf W}\cdot\overline{\psi}+{\bf V}\cdot\overline{\psi},
\end{eqnarray}

where ${\bf W}\in\Gamma\big(\mathbb{C}l_g(TU)\big)$ and ${\bf V}\in\Gamma(TU)$. With a little work, on can compute the expansion of ${\bf W}$ and ${\bf V}$ in a neighbourhood of $p\in U$. In fact, if $m\in U$ and $r$ denotes the distance from $m$ to $p$, we have:
\begin{eqnarray}
b^j_i & = & \delta_i^j-\frac{1}{6}R_{i\alpha\beta j}(p)x^\alpha x^\beta+O(r^3)\\
{\bf V} & = & \big(-\frac{1}{4}(Ric)_{\alpha k}(p)x^\alpha+O(r^2)\big)e_k\\
|{\bf W}| & =& O(r^3),
\end{eqnarray}

where $R_{ijkl}$ (resp. $(Ric)_{ik}$) are the components of the Riemann (resp. Ricci) curvature tensor. Now consider the smooth spinor field defined on $(V,\xi)$ by:
\begin{eqnarray*}
\psi(x)=f^{\frac{n}{2}}(x)(1-x)\cdot\psi_0
\end{eqnarray*}

where $f(x)=\frac{2}{1+r^2}$ (with $r^2=x_1^2+\cdots+x_n^2$) and $\psi_0\in\Sigma_\xi(V)$ is a constant spinor which can be chosen such that $|\psi_0|=1$. A straightforward computation shows that:
\begin{eqnarray}
D_\xi\psi=\frac{n}{2}f\psi,\quad|\psi|^2 =  f^{n-1}\quad{\text{and}}\quad |D_\xi\psi|^2=\frac{n^2}{4}f^{n+1}.
\end{eqnarray}

With these constructions, we can prove the main statement of this section.\\

{\it Proof of Theorem \ref{existance}:}
Let $\varepsilon>0$ and $\psi$ the spinor field described above, then we define:
\begin{eqnarray}
\psi_\varepsilon(x):=\eta\psi\big(\frac{x}{\varepsilon}\big)\in\Gamma\big(\Sigma_\xi(\mathbb{R}^n)\big)
\end{eqnarray}

where $\eta=0$ on $\mathbb{R}^n\setminus B_p(2\delta)$, $\eta=1$ on $B_p(\delta)$ and $0<\delta<1$ is chosen such that $B_p(2\delta)\subset V$. Since the support of the spinor field $\psi_\varepsilon$ lies in the open set $V$ of $\mathbb{R}^n$, one can use the trivialization described previously to obtain a spinor field $\overline{\psi_\varepsilon}$ over $(M,g)$. On the other hand, because of the conformal covariance of $\lambda_{\min}$, we can assume that the metric $g$ satisfies $Ric(p)_{ij}=0$. First we compute:
\begin{eqnarray*}
D_g\ov{\psi}_{\varepsilon}(x) & = & \ov{\nabla}\eta\cdot\ov{\psi}_{\varepsilon}\big(\frac{x}{\varepsilon}\big)+\frac{\eta}{\varepsilon}\frac{n}{2}f\big(\frac{x}{\varepsilon}\big)\ov{\psi}_{\varepsilon}\big(\frac{x}{\varepsilon}\big)+\eta\sum_{i,j}(b^j_i-\delta_i^j)\ov{\partial_i\cdot\nabla_{\partial_j}\big(\psi(\frac{x}{\varepsilon})\big)}\\
&  & + \eta{\bf W}\cdot\ov{\psi}_{\varepsilon}\big(\frac{x}{\varepsilon}\big)+\eta{\bf V}\cdot\ov{\psi}_{\varepsilon}\big(\frac{x}{\varepsilon}\big)
\end{eqnarray*}

where $|{\bf W}|=O(r^3)$ and $|{\bf V}|=O(r^2)$ (since $Ric(p)_{ij}=0$). Using \cite{amm3}, we have:
\begin{eqnarray*}
|D_g\ov{\psi}_\varepsilon|^2(x)\leq\frac{n^2}{4\varepsilon^2}f^{n+1}\big(\frac{x}{\varepsilon}\big)+Cr^4f^{n-1}\big(\frac{x}{\varepsilon}\big)+\frac{C}{\varepsilon}r^2f^{n-\frac{1}{2}}\big(\frac{x}{\varepsilon}\big)=\frac{n^2}{4\varepsilon^2}f^{n+1}\big(\frac{x}{\varepsilon}\big)\big(1+\Lambda(x)\big)
\end{eqnarray*}

where $\Lambda(x)=C\varepsilon^2r^4f^{-2}\big(\frac{x}{\varepsilon}\big)+C\varepsilon r^2f^{-\frac{3}{2}}\big(\frac{x}{\varepsilon}\big)$. Now note that for all $u\geq -1$:
\begin{eqnarray*}
(1+u)^{\frac{n}{n+1}}\leq 1+\frac{n}{n+1}u,
\end{eqnarray*}

then we get:
\begin{eqnarray}\label{devd}
|D_g\ov{\psi}_\varepsilon|^{\frac{2n}{n+1}}(x)\leq\Big(\frac{n}{2\varepsilon}\Big)^{\frac{2n}{n+1}}f^{n}\big(\frac{x}{\varepsilon}\big)+\frac{n}{n+1}\Big(\frac{n}{2\varepsilon}\Big)^{\frac{2n}{n+1}}f^{n}\big(\frac{x}{\varepsilon}\big)\Lambda(x).
\end{eqnarray}

On the other hand, since $p\in M$ is a point where $H$ is maximum, we have:
\begin{eqnarray*}
H(x)=H(p)+O(r^2)
\end{eqnarray*}

which yields:
\begin{eqnarray}\label{devf}
H(x)^{-\frac{n-1}{n+1}}=H(p)^{-\frac{n-1}{n+1}}\big(1+O(r^2)\big).
\end{eqnarray}

An integration combining (\ref{devd}) and (\ref{devf}) gives:
\begin{eqnarray*}
\int_{B_p(2\delta)}H^{-\frac{n-1}{n+1}}|D_g\ov{\psi}_\varepsilon|^{\frac{2n}{n+1}}dv(g)\leq\Big(\frac{n}{2\varepsilon}\Big)^{\frac{2n}{n+1}}H(p)^{-\frac{n-1}{n+1}}\big({\bf A}+{\bf B}+{\bf C}+{\bf D}\big)
\end{eqnarray*}

where:
\begin{eqnarray*}
{\bf A} & = & \int_{B_p(2\delta)}f^{n}\big(\frac{x}{\varepsilon}\big)dv(g)\\
{\bf B} & = & C\int_{B_p(2\delta)}f^{n}\big(\frac{x}{\varepsilon}\big)\Lambda(x)\,dv(g)\\
{\bf C} & = & C\int_{B_p(2\delta)}r^2f^{n}\big(\frac{x}{\varepsilon}\big)\,dv(g)\\
{\bf D} & = & C\int_{B_p(2\delta)}r^2f^{n}\big(\frac{x}{\varepsilon}\big)\Lambda(x)\, dv(g).
\end{eqnarray*}

Since the function $f$ is radially symmetric, we have:
\begin{eqnarray*}
{\bf A}  =  \int_{0}^{2\delta}f^{n}\big(\frac{x}{\varepsilon}\big)\omega_{n-1}G(r)r^{n-1}dr
\end{eqnarray*}

where: 
\begin{eqnarray*}
G(r)=\int_{\mathbb{S}^{n-1}}\sqrt{|g|_{rx}}d\sigma(x)\quad\text{with}\quad|g|_y:=\det g_{ij}(y).
\end{eqnarray*}

Now using the fact that $Ric_{ij}(p)=0$, one can compute that (see \cite{hebey}, for example):
\begin{eqnarray*}
G(r)\leq 1+O(r^4).
\end{eqnarray*}

Thus, a direct calculation shows that if $n\geq 3$:
\begin{eqnarray*}
{\bf A}=\omega_{n-1} I \varepsilon^n+o(\varepsilon^{n}),
\end{eqnarray*}

where $I=\int_0^{+\infty}r^{n-1}f^n(r)dr$. In the same way, we can prove that for $n\geq 3$:
\begin{eqnarray*}
{\bf B}={\bf C}={\bf D}=o(\varepsilon^{n}).
\end{eqnarray*}

In brief, we showed that:
\begin{eqnarray*}
\int_{B_p(2\delta)}H^{-\frac{n-1}{n+1}}|D_g\ov{\psi}_\varepsilon|^{\frac{2n}{n+1}}dv(g)=\big(\frac{n}{2}\big)^{\frac{2n}{n+1}}\big(\omega_{n-1} I\big) H(p)^{-\frac{n-1}{n+1}}\varepsilon^{\frac{n(n-1)}{n+1}}\big(1+o(1)\big),
\end{eqnarray*}

hence:
\begin{eqnarray*}
\Big(\int_M H^{-\frac{n-1}{n+1}}|D_g\ov{\psi}_\varepsilon|^{\frac{2n}{n+1}}dv(g)\Big)^{\frac{n+1}{n}}=\big(\frac{n}{2}\big)^{2}\big(\omega_{n-1} I\big)^{\frac{n+1}{n}} H(p)^{-\frac{n-1}{n}}\varepsilon^{n-1}\big(1+o(1)\big).
\end{eqnarray*}

The denominator of the functional $\lambda_{\min}$ can also be estimated and similar computations give (see \cite{amm3}): 
\begin{eqnarray*}
\int_M\<D_g\ov{\psi}_\varepsilon,\ov{\psi}_\varepsilon\>dv(g)=\frac{n}{2}\omega_{n-1} I \varepsilon^{n-1}+o(\varepsilon^{n+1})
\end{eqnarray*}

for $n\geq 3$. Combining these two estimates yields:
\begin{eqnarray*}
\lambda_{\min}\leq K(n)^{-1}\big(\max_M H\big)^{-\frac{2}{p_D}}\big(1+o(1)\big)
\end{eqnarray*}

which concludes the proof.
\hfill$\square$\\

\begin{remark}
We can derive a similar result for the case of $2$-dimensional manifolds. It is sufficient to adapt the proof of \cite{amm3} in our situation and one can show that if $(M^2,g)$ is a smooth surface we have:
\begin{eqnarray*}
\lambda_{\min}\leq 2\sqrt{\pi}\,(\max_M H)^{-2}.
\end{eqnarray*}
\end{remark}

\begin{remark}
This result is in the spirit of the one obtained by Aubin in \cite{aubin:76} for the conformal Laplacian. Indeed, in this article, the author proves that on an $n$-dimensional compact Riemannian manifold with $n>4$, if $f,h$ are smooth positive functions on $M$ such that:
\begin{eqnarray*}
h(p)-R_g(p)+\frac{n-4}{2}\frac{\Delta_g f(p)}{f(p)}<0,
\end{eqnarray*}

\noindent where $f(p)=\max_{x\in M}f(x)$, then the nonlinear equation:
\begin{eqnarray*}
4\frac{n-1}{n-2}\Delta_g u+h u= f u^{\frac{n+2}{n-2}}
\end{eqnarray*}

\noindent admits a smooth positive solution. We could hope to obtain a similar criterion for the equation studied in this paper. However, if one carries out the computations in the proof of Theorem \ref{existance}, we obtain for $n\geq 5$:
\begin{eqnarray*}
\lambda_{\min}= K(n)^{-1}\big(\max_M H\big)^{-\frac{2}{p_D}}\big(1+\frac{n-1}{2n(n-2)}\frac{\Delta H(p)}{H(p)}\,\varepsilon^{2}+o(\varepsilon^{2})\big).
\end{eqnarray*}

Thus no conclusion could be made since at a point $p\in M$ where $H$ is maximum we have $\Delta H(p)\geq 0$.
\end{remark}

%%%%%%%%%%%%%%%%%%%%%%%%%%%%%%%%%%%%%%%%%%%%%%%%%%%%%%%%%%%%%%%%%%%%%%%%%%%%%%%%%%%%%%%%%%%%%%%%%%%%%%%%

\subsection{An existence result}\label{parthree}

%%%%%%%%%%%%%%%%%%%%%%%%%%%%%%%%%%%%%%%%%%%%%%%%%%%%%%%%%%%%%%%%%%%%%%%%%%%%%%%%%%%%%%%%%%%%%%%%%%%%%%%%

To end this section, we give conditions on the manifold $(M^n,g)$ and on the function $H\in C^\infty(M)$ which ensure that (\ref{strictdirac}) is fulfilled. Then applying Theorem~\ref{existence}, we get the existence of a solution to the nonlinear Dirac equation (\ref{nld}). The condition on $H$ is a technical one given by:
\begin{equation}\label{conditionH}
\begin{array}{ll}
\textrm{ {\it There is a maximum point }} p\in M \;\textrm{{\it at which all partial derivaties of }} H\\
 \textrm{{\it of order less than or equal to }} (n-1) \;\text{{\it vanish.}}
\end{array}
\end{equation}

The result we obtain is the following:
\begin{theorem}\label{EXISTANCE}
Let $(M^n,g)$ be an $n$-dimensional compact Riemannian spin manifold. Assume that $(M^n,g)$ is locally conformally flat and $H\in C^\infty(M)$ a smooth positive function on $M$ for which (\ref{conditionH}) holds. Then if the Dirac operator is invertible and the mass endomorphism has a positive eigenvalue, there exists a spinor field solution of the nonlinear Dirac equation (\ref{nld}).
\end{theorem}

This result is quite close to the work of Escobar and Schoen \cite{escobarschoen} and relies on the construction of Ammann, Humbert and Morel \cite{amm4} of the mass endomorphism.

We first briefly recall the construction of the mass endomorphism. For more details, we refer to \cite{amm4}. Consider a point $p\in M$ and suppose that there is a neighborhood $U$ of $p$ which is flat. Since we assumed that the Dirac operator has a trivial kernel, one can show that the Green function $G_D$ of the Dirac operator has the following expansion in $U$:
\begin{eqnarray*}
\omega_{n-1}G_D(x,p)\psi_0=-\frac{x-p}{|x-p|^n}\cdot\psi_p+v(x,p)\psi_p
\end{eqnarray*}

for all $x\in U$ and where $v(\,.\,,p)\psi_p$ is a smooth harmonic spinor near $p$ with $\psi_p\in\Sigma_p(M)$. The mass endomorphism is then the self-adjoint endomorphism of the fiber $\Sigma_p(M)$ defined by:
\begin{eqnarray*}
\alpha_p(\psi_p)=v(p,p)\psi_p.
\end{eqnarray*}

This operator shares many properties with the mass of the Green function of the conformal Laplacian. One of them is that the sign of its eigenvalues is invariant under conformal changes of metrics which preserves the flatness near $p$. With this construction, we can prove the main result of this section.\\

{\it Proof of Theorem \ref{EXISTANCE}:}
We have to construct a test-spinor which will be estimated in the variational characterization of $\lambda_{\min}$. The assumption on the mass endomorphism implies that (\ref{strictdirac}) is fulfilled and the result will follow from Theorem \ref{existence}. The test-spinor is exactly the one used in \cite{amm4}. In order to make this paper self-contained, we have chosen to briefly recall this construction. First, since $\lambda_{\min}$ is a conformal invariant of $(M^n,g)$ which is locally conformally flat, one can suppose that the metric is flat near a point $p\in M$ where (\ref{conditionH}) is satisfied. Now for $\varepsilon>0$ we set:
\begin{eqnarray*}
\xi:=\varepsilon^{\frac{1}{n+1}}\qquad\varepsilon_0:=\frac{\xi^n}{\varepsilon}f\big(\frac{\xi}{\varepsilon})^{\frac{n}{2}}
\end{eqnarray*}

where $f(r)=\frac{2}{1+r^2}$ is the function used in the previous section. The test-spinor is then defined by:
$$\Phi_{\varepsilon}(x)=
\left\lbrace
\begin{array}{ll}
f\big(\frac{x}{\varepsilon}\big)^{\frac{n}{2}}\big(1-\frac{x}{\varepsilon}\big)\cdot\psi_p+\varepsilon_0\alpha_p(\psi_p) \quad & \quad \textrm{if }r\leq\xi\\ \\
\varepsilon_0\big(\omega_{n-1}G_D(x,p)-\eta(x)\theta_p(x)\big)+\eta(x) f\big(\frac{\xi}{\varepsilon}\big)^{\frac{n}{2}}\psi_p\quad & \quad\textrm{if } \xi\leq r\leq 2\xi\\ \\
\varepsilon_0\,\omega_{n-1}G_D(x,p) \quad & \quad\textrm{if } r\geq 2\xi
\end{array}
\right.$$

where $\eta$ is a cut-off function such that:
$$\eta=\left\lbrace\begin{array}{ll}
1 & \quad\text{on}\;B_p(\xi)\\
0 & \quad\text{on}\; M\setminus B_p(2\xi)
\end{array}
\qquad\text{and}\qquad|\nabla\eta|\leq\frac{2}{\xi}
\right.$$

and $\theta_p(x):=v(x,p)\psi_p-\alpha_p(\psi_p)$ is a smooth spinor field (harmonic near $p$) which satisfies $|\theta_p|=O(r)$. Now an easy calculation shows that:
$$|D\Phi_{\varepsilon}|^{\frac{2n}{n+1}}(x)=
\left\lbrace
\begin{array}{ll}
\big(\frac{n}{2}\big)^{\frac{2n}{n+1}}\varepsilon^{-\frac{2n}{n+1}}f\big(\frac{r}{\varepsilon}\big)^n \quad & \quad \textrm{if }r\leq\xi\\ \\
|\varepsilon_0\nabla\eta(x)\cdot\theta_p(x)-f\big(\frac{\xi}{\varepsilon}\big)^{\frac{n}{2}}\nabla\eta(x)\cdot\psi_p|^{\frac{2n}{n+1}}\quad & \quad\textrm{if } \xi\leq r\leq 2\xi\\ \\
0 \quad & \quad\textrm{if } r\geq 2\xi.
\end{array}
\right.$$

On the other hand, since the function $H$ satisfies the condition (\ref{conditionH}), we get:
\begin{eqnarray*}
H(x)=H(p)+O(r^n)
\end{eqnarray*}

that is:
\begin{eqnarray*}
H(x)^{-\frac{n-1}{n+1}}=H(p)^{-\frac{n-1}{n+1}}\big(1+O(r^n)\big).
\end{eqnarray*}

We can now give the estimate of the functional (\ref{lambdaq}) (with $q=q_D$) evaluated at the spinor field $\Phi_\varepsilon$. First, on $B_p(\xi)$ we have:
\begin{eqnarray*}
\int_{B_p(\xi)}H^{-\frac{n-1}{n+1}}|D\Phi_\varepsilon|^{\frac{2n}{n+1}}dx\leq\big(\frac{n}{2}\big)^{\frac{2n}{n+1}}\varepsilon^{-\frac{2n}{n+1}}H(p)^{-\frac{n-1}{n+1}}\Big(\int_{B_p(\xi)}f^n\big(\frac{x}{\varepsilon}\big)dx+C\int_{B_p(\xi)}r^nf^n\big(\frac{x}{\varepsilon}\big)dx\Big)
\end{eqnarray*} 

and we compute that:
\begin{eqnarray*}
\int_{B_p(\xi)}f^n\big(\frac{x}{\varepsilon}\big)dx & \leq & \varepsilon^n\int_{\mathbb{R}^n}f^n(x)dx\\
\int_{B_p(\xi)}r^nf^n\big(\frac{x}{\varepsilon}\big)dx & = & o(\varepsilon^{2n-1}).
\end{eqnarray*}

Finally we obtain:
\begin{eqnarray*}
\int_{B_p(\xi)}H^{-\frac{n-1}{n+1}}|D\Phi_\varepsilon|^{\frac{2n}{n+1}}dx & = & \big(\frac{n}{2}\big)^{\frac{2n}{n+1}}\varepsilon^{\frac{n(n-1)}{n+1}}H(p)^{-\frac{n-1}{n+1}}I\,\big(1+o(\varepsilon^{n-1})\big)
\end{eqnarray*}

where $I=\int_{\mathbb{R}^n}f^n(x)dx$. On $C_p(\xi):=B_p(2\xi)\setminus B_p(\xi)$:
\begin{eqnarray*}
\int_{C_p(\xi)}H^{-\frac{n-1}{n+1}}|D\Phi_\varepsilon|^{\frac{2n}{n+1}}dx & \leq & C\int_{C_p(\xi)}|\varepsilon_0\nabla\eta\cdot\theta_p|^{\frac{2n}{n+1}}dx+C\int_{C_p(\xi)}|f\big(\frac{\xi}{\varepsilon}\big)^{\frac{n}{2}}\nabla\eta\cdot\psi_p|^{\frac{2n}{n+1}}dx\\
& + & C\int_{C_p(\xi)}r^n|\varepsilon_0\nabla\eta\cdot\theta_p|^{\frac{2n}{n+1}}dx
+C\int_{C_p(\xi)}r^n|f\big(\frac{\xi}{\varepsilon}\big)^{\frac{n}{2}}\nabla\eta\cdot\psi_p|^{\frac{2n}{n+1}}dx
\end{eqnarray*}

and since $\varepsilon_0\leq C\varepsilon^{n-1}$, $|\nabla\eta|\leq 2\xi^{-1}$, $|\theta_p|=O(r)$ and $Vol\big(C_p(\xi)\big)\leq C\xi^n$, we get:
\begin{eqnarray*}
\int_{C_p(\xi)}H^{-\frac{n-1}{n+1}}|D\Phi_\varepsilon|^{\frac{2n}{n+1}}dx =o\big(\varepsilon^{\frac{(2n+1)(n-1)}{n+1}}\big).
\end{eqnarray*}

In conclusion, the numerator of (\ref{lambdaq}) is given by:
\begin{eqnarray*}
\Big(\int_{M}H^{-\frac{n-1}{n+1}}|D\Phi_\varepsilon|^{\frac{2n}{n+1}}dv(g)\Big)^{\frac{n+1}{n}}= \big(\frac{n}{2}\big)^{2}\varepsilon^{n-1}H(p)^{-\frac{n-1}{n}}
I^{\frac{n+1}{n}}\,\big(1+o(\varepsilon^{n-1})\big)
\end{eqnarray*}

Similar computations for the denominator lead to (see also \cite{amm4}):
\begin{eqnarray*}
\int_M\<D\Phi_{\varepsilon},\Phi_{\varepsilon}\>dv(g)=\frac{n}{2}\varepsilon^{n-1}I\big(1+J\,\<\psi_p,\alpha_p(\psi_p)\>\varepsilon^{n-1}+o(\varepsilon^{n-1})\big)
\end{eqnarray*}

where $J=\int_{\mathbb{R}^n}f(x)^{\frac{n}{2}+1}dx$. Now we choose $\psi_p\in\Sigma_p(M)$ as an eigenspinor for the mass endomorphism associated with a positive eigenvalue $\lambda$ and we finally get:
\begin{eqnarray*}
\lambda_{\min}\leq \mathcal{F}_{q_D}(\Phi_\varepsilon)=K(n)^{-1}\big(\max_M H\big)^{-\frac{n-1}{n}}\big(1-\lambda J\,\varepsilon^{n-1}+o(\varepsilon^{n-1})\big).
\end{eqnarray*}

Now it is clear that for $\varepsilon>0$ sufficiently small, (\ref{strictdirac}) is true and thus Theorem \ref{existence} allows to conclude.
\hfill$\square$

\begin{remark}
In dimension two, a Riemannian surface is always locally conformally flat and condition (\ref{conditionH}) is satisfied for all $H\in C^{\infty}(M)$ however the mass endomorphism vanishes (see \cite{amm4}) and so Theorem~\ref{EXISTANCE} cannot be applied.
\end{remark}

%%%%%%%%%%%%%%%%%%%%%%%%%%%%%%%%%%%%%%%%%%%%%%%%%%%%%%%%%%%%%%%%%%%%%%%%%%%%%%%%%%%%

\section{A remark on manifolds with boundary}

%%%%%%%%%%%%%%%%%%%%%%%%%%%%%%%%%%%%%%%%%%%%%%%%%%%%%%%%%%%%%%%%%%%%%%%%%%%%%%%%%%%%

In this last section, we briefly study the case of manifolds with boundary. Since the calculations are quite close to those of the boundaryless case, we only point out arguments which need some explanations.  

Indeed, let $(M^n,g)$ be an $n$-dimensional compact Riemannian spin manifold with smooth boundary equipped with a chirality operator $\gamma$, that is an endomorphism of the spinor bundle which satisfies:
\begin{eqnarray*}\label{chi}
\gamma^2 =Id, & \quad\<\gamma\psi,\gamma\varphi\>=\<\psi,\varphi\>\\
\nabla_{X}(\gamma\psi)=\gamma(\nabla_{X}\psi), & \quad X\cdot\gamma\psi=-\gamma(X\cdot\psi),
\end{eqnarray*}

\noindent for all $X\in\Gamma(T M)$ and for all spinor fields $\psi,\varphi\in\Gamma\big(\Sigma_g(M)\big)$. The orthogonal projection:
\begin{eqnarray*}
B^\pm_g:=\frac{1}{2}(Id\pm\nu_g\cdot\gamma),
\end{eqnarray*}

where $\nu_g$ denotes the inner unit vector fields normal to $\partial M$, defines a (local) elliptic boundary condition (called the chiral bag boundary condition or $(CHI)$ boundary condition) for the Dirac operator $D_g$ of $(M,g)$. Moreover, under this boundary condition, the spectrum of the Dirac operator consists of entirely isolated real eigenvalues with finite multiplicity. In \cite{sr3} (see also \cite{mathese}), we define a spin conformal invariant similar to (\ref{sci}) using this boundary condition. More precisely, if $\lambda^{\pm}_1(\ov{g})$ stands for the first eigenvalue of the Dirac operator $D_{\ov{g}}$ under the chiral bag boundary condition $B^\pm_{\ov{g}}$ then the chiral bag invariant is defined by:
\begin{eqnarray*}
\lambda_{\min}(M,\pa M):=\underset{\ov{g}\in[g]}{\inf}|\lambda^{\pm}_1(\ov{g})|\Vol(M,\ov{g})^{\frac{1}{n}}
\end{eqnarray*}

and one can check that:
\begin{eqnarray}\label{chirbagcar}
\lambda_{\min}(M,\pa M)=\underset{\varphi\neq 0}{\inf}\Big\{\frac{\big(\int_{M}|D_g\varphi|^{\frac{2n}{n+1}}dv(g)\big)^{\frac{n+1}{n}}} {\big|\int_{M}\<D_g\varphi,\varphi\>dv(g)\big|}\Big\}
\end{eqnarray}

where the infimum is taken for all spinor fields $\varphi\in H_1^{q_D}$ such that $B^\pm_g\varphi_{|\partial M}=0$. On the round hemisphere $(\mathbb{S}^n_+,g_{\rm{st}})$, we can compute that:
\begin{eqnarray}\label{ahm}
\lambda_{\min}(\mathbb{S}^n_+,\pa\hs)=\frac{n}{2}\Big(\frac{\omega_n}{2}\Big)^{\frac{1}{n}}=2^{-(1/n)}{K}(n)^{-1},
\end{eqnarray}
 
and using the conformal covariance of (\ref{chirbagcar}) and the fact that the hemisphere is conformally isometric to the half Euclidean space $(\mathbb{R}^n_+,\xi)$, we conclude that:
\begin{eqnarray}\label{inreb}
\Big|\int_{\mathbb{R}^n_+}\<D_\xi\psi,\psi\>dx\Big|\leq {2^{\frac{1}{n}}K(n)}\Big(\int_{\mathbb{R}^n_+}|D_\xi\psi|^{\frac{2n}{n+1}}dx\Big)^{\frac{n+1}{n}}
\end{eqnarray}

for all $\psi\in\Gamma_c\big(\Sigma_{\xi}(\mathbb{R}^n_+)\big)$ where $\Gamma_c\big(\Sigma_{\xi}(\mathbb{R}^n_+)\big)$ denotes the space of smooth spinor fields over $(\mathbb{R}^n_+,\xi)$ with compact support. In order to prove a Sobolev-type inequality for manifolds with boundary, we give a result similar to Lemma \ref{pdo} in this context:
\begin{lemma}\label{pdob}
If the Dirac operator is invertible under the chiral bag boundary condition then there exists a constant $C>0$ such that:
\begin{eqnarray*}
||\varphi||_{p_D}\leq C ||D_g\varphi||_{q_D}
\end{eqnarray*}

for all $\varphi\in H^{q_D}_1$ such that $B^\pm_{g}\varphi_{|\pa\M}=0$.
\end{lemma}

{\it Proof:} 
Since the Dirac operator is assumed to be invertible and since the Fredholm property of $D_g$ does not depend on the choice of the Sobolev spaces (see \cite{schwartz}), we have that:
\begin{eqnarray*}
D_g:\mathcal{H}^{q_D}_\pm:=\big\{\varphi\in H^{q_D}_1\,/\,B^\pm_{g}\varphi_{|\pa\M}=0\big\}\longrightarrow L^{q_D}
\end{eqnarray*}

defines a continuous bijection. Using the open mapping theorem, the inverse map is also continuous and then we get the existence of a constant $C>0$ such that:
\begin{eqnarray*}
||\varphi||_{H^{q_D}_1}= ||D^{-1}_g(D_g\varphi)||_{H^{q_D}_1}\leq C||D_g\varphi||_{q_D},
\end{eqnarray*}

for all $\varphi\in\mathcal{H}^{q_D}_\pm$. On the other hand, the Sobolev embedding Theorem implies that the map $H^{q_D}_1\hookrightarrow L^{p_D}$ is continuous, so there exists a constant $C>0$ such that:
\begin{eqnarray*}
||\varphi||_{p_D}\leq C||\varphi||_{H^{q_D}_1},
\end{eqnarray*}

for all $\varphi\in H_1^{q_D}$ and this concludes the proof.
\hfill$\square$\\

\begin{remark}\label{re2}
Lemma~\ref{pdob} gives a result similar to Lott's one (see \cite{lott}) for the Dirac operator on manifolds with boundary, that is if $D_g$ is invertible under the chiral bag boundary condition, then:
\begin{eqnarray}\label{lottbord}
\lambda_{\min}(M,\pa M)>0.
\end{eqnarray}
 
\noindent Indeed, the H\"older inequality gives:
\begin{eqnarray*}
\Big|\int_{M}\<D_g\varphi,\varphi\>dv(g)\Big|\leq ||\varphi||_{p_D}||D_g\varphi||_{q_D},
\end{eqnarray*}

\noindent and then Lemma~\ref{pdo} yields:
\begin{eqnarray*}
\frac{||D_g\varphi||_{q_D}^2}{\Big|\int_{M}\<D_g\varphi,\varphi\>dv(g)\Big|}\geq \frac{||D\varphi||_{q_D}}{||\varphi||_{p_D}}\geq C,
\end{eqnarray*}

\noindent for all $\varphi\in\mathcal{H}^{q_D}_\pm$. Using the variational characterization (\ref{chirbagcar}) of $\lambda_{\min}(M,\pa M)$ leads to the result. In \cite{sr}, we give an explicit lower bound for the chiral bag invariant given by:
\begin{eqnarray}\label{hijbor}
\lambda_{\mathrm{min}}(M,\pa M)^2\geq\frac{n}{4(n-1)}\mu_{[g]}(M,\pa M).
\end{eqnarray}

The number $\mu_{[g]}(M,\pa M)$ is a conformal invariant of the manifold introduced by Escobar in \cite{escobar:92} to study the Yamabe problem on manifolds with boundary and defined by:
\begin{eqnarray*}
\mu_{[g]}(M,\pa M)=\underset{u\in C^{1}(M),u\neq 0}{\mathrm{inf}}\frac{\int_{\M}\big(4\frac{n-1}{n-2}|\na u|^2+R_g u^2\big)dv(g)+2(n-1)\int_{\pa M}h_g u^2 ds(g)}{\Big(\int_{M}u^{N} ds(g)\Big)^{\frac{2}{N}}}.
\end{eqnarray*}

This invariant is called the Yamabe invariant of $(M^n,g)$. Here $h_g$ denotes the mean curvature of the boundary of $(\partial M,g)$ in $(M,g)$. Inequality~(\ref{hijbor}) is significant only if the Yamabe invariant is positive and in this case, the Dirac operator under the chiral bag boundary condition is invertible. So Inequality (\ref{lottbord}) is more general than (\ref{hijbor}) however it does not give an explicit lower bound.
\end{remark}

We can now argue like in the proof of Theorem~\ref{sli} and state a Sobolev-like inequality on manifolds with boundary:
\begin{theorem}\label{slib}
Let $(M^n,g,\sigma)$ be an $n$-dimensional compact spin manifold with a non empty smooth boundary and equipped with a chirality operator. Moreover, we assume that the Dirac operator under the chiral bag boundary condition is invertible. Then for all $\varepsilon>0$, there exists a constant $B_{\varepsilon}$ such that:
\begin{eqnarray*}
\Big|\int_{M}\<D_g\varphi,\varphi\>dv(g)\Big|\leq \Big(2^{1/n}{K}(n)+\varepsilon\Big)\Big(\int_{M}|D_g\varphi|^{q_D}dv(g)\Big)^{2/q_D}+B_{\varepsilon}\Big(\int_{M}|\varphi|^{q_D}dv(g)\Big)^{2/q_D},
\end{eqnarray*}

for all $\varphi\in H^{q_D}_1$ such that $B^\pm_{g}\varphi_{|\pa\M}=0$.
\end{theorem}

%%%%%%%%%%%%%%%%%%%%%%%%%%%%%%%%%%%%%%%%%%%%%%%%%%%%%%%%%%%%%%%%%%%%%%%%%%%%%%%%%%%%

\bibliographystyle{amsalpha}     
\bibliography{biblio}

%%%%%%%%%%%%%%%%%%%%%%%%%%%%%%%%%%%%%%%%%%%%%%%%%%%%%%%%%%%%%%%%%%%%%%%%%%%%%%%%%%%%

%%%%%%%%%%%%%%%%%%%%%%%%%%%%%%%%%%%%%%%%%%%%%%%%%%%%%%%%%%%%%%%%%%%%%%%%%%%%%%%%%%%%     

\end{document}